\documentclass[14pt]{extarticle}
\usepackage[margin=1in]{geometry}
\usepackage{dsfont}
\usepackage{amsfonts}
\usepackage{tikz}
\usetikzlibrary{matrix}

\newtheorem{theorem}{Theorem}[section]
\newtheorem{pro}{Proposition}[section]
\newtheorem{lemma}{Lemma}[section]

\newtheorem{cor}{Corollary}[section]

\newcommand{\proof}[1]{\noindent{\it\bf Proof:#1\ }}
\newcommand{\QED}{\hfill$\Box$\medskip}

\begin{document}

\title{ Geometric Equivariant Extension of Sections 
	  in GW Theory I}
\author{
Gang Liu }
\date{June 22,  2015}
\maketitle

\section{Introduction}

Let ${\widetilde {\cal B}}=:{\widetilde {\cal B}}_{k, p}$  be the  space of stable $L_k^p$-maps $f:S^2\rightarrow M$ from the Riemann sphere to a compact  almost 
	complex manifold $(M, J)$. Consider the bundle  ${\widetilde {\cal L}}_{k-1, p}\rightarrow {\widetilde {\cal B}_{k, p}} $ with the fiber $({\widetilde {\cal L}}_{k-1, p})_f=L_{k-1}^p(S^2, \wedge^{0, 1}_J(f^*TM)) $.
	 It is well-known that $ ({\widetilde {\cal L}}_{k-1, p},  {\widetilde {\cal B}_{k, p}})  $ is a Banach bundle of class $C^{\infty}.$ The reparametrization group $G=:{\bf PSL}(2, {\bf C})$ acting on 
	$ ({\widetilde {\cal L}}_{k-1, p},  {\widetilde {\cal B}_{k, p}})  $ continuously with local slices $S_f$ for any $f\in {\widetilde {\cal B}_{k, p}}$. Since the $G$-action  is only continuous, the coordinate transformations between these local slices $S_f$ as well as the transition functions between the local bundles $ ({\widetilde {\cal L}}_{k-1, p}|{S_f}, S_f)  $ are only continuous in general. Because of this well-known  difficulty in  lack of differentiability, to establish the analytic foundation 
	of GW theory by the method of [L] it is crucial to have  sufficiently many smooth sections of the local bundle
	${\cal L}_{k-1, p}|_{S_f}\rightarrow S_f$  that are  still smooth viewed in any other slices $S_{h}$.
	
	In [L], we have  given two different methods to construct such sections starting from an element $\xi\in C^{\infty}(S^2, \wedge^{0, 1}(f^*TM))
	\subset L_{k-1}^p(S^2, \wedge^{0, 1}(f^*TM))=({\widetilde {\cal L}}_{k-1, p})_f$ for  a smooth  map $f$. The purpose  of this paper is to explain one of the constructions,  the geometric 
	$G$-equivariant extension ${\xi}_{{\cal O}_{S_f}}$ of $\xi$, and  to improve
	 the results in [L] on this construction.
	Throughout this paper, we will assume that   $k-2/p>1$. Note that under this condition  the $L_{k-1}^p$ functions on $S^2$ form a Banach algebra. Set $m_0=[k-2/p].$
	
	First observe that for a fixed $\phi\in G$, the action on $ ({\widetilde {\cal L}}_{k-1, p},  {\widetilde {\cal B}_{k, p}})  $  is a $C^{\infty}$
	automorphism. Hence for the purpose of this paper and its sequels, we can replace $G$ by $G_e$ and only consider the local $G_e$-action. Here    $G_e$  is  a local chart of $G$ containing the identity $e\in G$ chosen  as follows.
Let ${\widetilde W}(f) $ be a local coordinate chart of ${\widetilde {\cal B}_{k, p}}$ centered at a smooth stable map $f$ such that the bundle ${\widetilde {\cal L}}_{k-1, p}|_{{\widetilde W}(f)}\rightarrow {\widetilde W}(f)$ is trivialized. Let $\Gamma_f$ be the isotropy group of $f$. Since $\Gamma_f$ is  finite, we can choose $G_e$ such that $G_e\cap \Gamma_f=\{e\}.$
Then   the $G_e$-action on ${\widetilde W}(f)$ is free for sufficiently small ${\widetilde W}(f)$. Let $S_f\subset {\widetilde W}(f)$ be a local slice of the action.  By shrinking $ {\widetilde W}(f)$ if it is necessary, 	we may  assume further that  the local $G_e$-orbit ${\cal O}_{S_f}$ of $S_f$ is equal to $ {\widetilde W}(f).$

	We now recall the two  constructions of the $G$-equivariant extensions  of $\xi$ in  [L].
	
	The first one is obtained as follows.  Under the above given trivialization, the given element $\xi$ in the central fiber can be considered as  a "constant", hence smooth section over  
	${\widetilde W}(f) $. Denote its restriction to the local slice $S_f$ by $\xi_{S_f}$ as a smooth section of ${\cal L}_{k-1, p}|_{S_f}\rightarrow S_f$.  
	Then using the local $G_e$-action, we  obtain a $G_e$-equivariant section ${\xi}_{{\cal O}_{S_f}}$ over $ {\widetilde W}(f)={\cal O}_{S_f}.$ More specifically, $\xi_{{\cal O}_{S_f}}$ is defined by 
	${\xi}_{{\cal O}_{S_f}}(h)=T(h)^*(\xi_{S_f}(h\circ T(h))).$ Here the map  $T$ is defined in the following proposition proved in next section.
	
	\begin{pro}
 Given  a smooth stable map $f:S^2\rightarrow M$, let   $S_f$ be the local slice $S_f^E$ or $S_f^{L_2}$ defined by the evaluation map or by the $L_2$ orthogonal complement to the orbit ${\cal O}_f$ of $f$. Then there is a $C^{m_0}$ or $C^{\infty}$ map $T: {\widetilde W}(f)\rightarrow G_e$ accordingly such that for any  $h\in {\widetilde W}(f)$, $h\circ T(h)\in S_f.$ Here  $m_0= [k-2/p]$. 	
	\end{pro}
	
	It was proved by a direct computation in [L] that  
	$ \xi_{{\cal O}_{S_f}}$ is of class $C^1$. A more conceptual
 proof of this will be given in the forth coming paper [L?], in which we will also show that generically ${\xi}_{{\cal O}_{S_f}}$ is exact of class $C^{1}$.
In contrast, we will show in this paper that the extension ${\xi}_{{\cal O}_{S_f}}$ by the 
second construction is of class $C^{m_0}$ or $C^{\infty}$. In other words, it has  the same degree of the smoothness as $T$ has.

To describe the second construction,  let $\cup_{i=1}^l V_i=M$ be a fixed open covering of $M$ . Consider a  corresponding   covering  $\cup_{i=1}^l D_i=S^2$ of $S^2$   satisfying  the condition that   $h(D_i) \subset V_i, i=1, \cdots, l$ for any $h\in  {{\widetilde { W}}(f)}$. In particular, 
for $(h, \phi)\in S_f\times G_e$, $h\circ \phi (D_i) \subset V_i.$
Fix a partition of unit $\{\alpha_i, i=1, \cdots, l\} $ on $S^2$ subordinate to the covering $\{D_i, i=1, \cdots, l\}.$ 
   Let $\beta_i, i=1, \cdots, l$ be  the cut-off functions defined on $V'_i$   such that 
   $\beta_i=1$ on $V_i$, where $V_i\subset \subset V_i'.$   %For each $i$, let $\alpha_i'$ be a cut-off function on $S^2$ supported on $D'_i$  with $D_i\subset \subset D_i'$ and $\alpha'_i=1$ on $D_i$.
Fix a smooth local  ${\bf C}$-frame ${\bf t}_{i}=\{t_{i1}, \cdots, t_{im}\}$ of $(TM, J)$ on
$V_i'.$  Then for any $\xi\in C^{\infty }(S^2, \wedge^{0, 1}(f^* TM))$ containing in the central fiber ${\widetilde {\cal L} }_f$, 

 $$\xi=\Sigma_i\alpha_i\cdot \xi=\Sigma_{i, \nu}\gamma_i^{\nu}\cdot t_{i\nu}\circ
f =\Sigma_{i, \nu}\gamma_i^{\nu}\cdot (\beta_i t_{i\nu}\circ
f).$$  Here $\gamma_i^{\nu}\in C^{\infty }(S^2, \wedge^{0, 1}_{S^2})$ supported on $D_i$ with $\alpha_i\cdot \xi=\Sigma_{ \nu}\gamma_i^{\nu}\cdot t_{i\nu}\circ
f $,  and $\beta_i t_{i\nu}\circ
f\in C^{\infty }(S^2, f^*TM)$ supported on $f^{-1}(V'_i)$. 

We require that the $G$-equivariant extension $\xi_{{\cal O}_{S_f}}$ to be defined satisfying  the condition that  $$\xi_{{\cal O}_{S_f}}=\Sigma_i(\alpha_i\cdot \xi)_{{\cal O}_{S_f}}
=\Sigma_{i, \nu}(\gamma_i^{\nu})_{{\cal O}_{S_f}}\cdot (\beta_i t_{i\nu}\circ
f)_{{\cal O}_{S_f}}.$$ 

Thus we only need to define the $G$-equivariant extensions $(\gamma)_{{\cal O}_{S_f}}$  for $\gamma \in C^{\infty }(S^2, \wedge^{0, 1}_{S^2})$ and $(\beta_i t_{i\nu}\circ
f)_{{\cal O}_{S_f}}.$

To get these  extensions, note that in above we have used the decomposition
$L_{k-1}^p(S^2, \wedge^{0, 1}(f^* TM))=L_{k-1}^p(S^2, \wedge^{0, 1})\otimes 
L_{k-1}^p(S^2, f^* TM).$  Here  we have used the fact that $L_{k-1}^p(S^2, \wedge^{0, 1})$ and $  
L_{k-1}^p(S^2, f^* TM)$ are the modules over the  complex Banach  algebra $L_{k-1}^p(S^2, {\bf C}), $ and tensor product is taken over  $L_{k-1}^p(S^2, {\bf C}).$
 A family version of this isomorphism gives rise the isomorphism of the bundles ${\widetilde {\cal L}}_{k-1, p}\simeq {\bf  \Omega}^{0, 1}_{k-1, p}\otimes {\widetilde {\cal T}}_{k-1, p}$. Here 
${\bf  \Omega}^{0, 1}_{k-1, p}= {\widetilde {\cal B}}_{k, p}\times L_{k-1}^p(S^2, \wedge^{0, 1})\rightarrow {\widetilde {\cal B}}_{k, p}$ is the trivial bundle and   ${\widetilde {\cal T}}_{k-1, p}$ is the tangent bundle 
$T{\cal B}_{k-1, p}$ restricted to  $ {\cal B}_{k, p}$ but with the "standard" bundle structure obtained by using   the  $J$-invariant  parallel transport on $TM$.  One can show (in next section) that this bundle structure is $C^{\infty}$ equivalent to that of $T{\cal B}_{k-1, p}$. However, the equivalence is not with respect to the induced
complex structure by $J$ on $T{\cal B}_{k-1, p}$. 

Using the trivialization of $\Omega^{0, 1}_{k-1, p}$ above, we get the constant
 extension $\gamma_{S_f}$ as a section of the bundle $\Omega^{0, 1}_{k-1, p}|_{S_f}\rightarrow S_f$. 
Then  $\gamma_{{\cal O}_{S_f}}: {\widetilde W}(f)\rightarrow {\bf  \Omega}^{0, 1}_{k-1, p}|_{{\widetilde W}(f)}=:{\widetilde W}(f)\times L_{k-1}^p(S^2, \wedge^{0, 1})$ is defined to be 
$\gamma_{{\cal O}_{S_f}}(h)=T(h)^*(\gamma_{S_f}(h\circ T(h)))$. Let $[\gamma_{S_f}]:S_f\rightarrow L_{k-1}^p(S^2, \wedge^{0, 1})$ and $[\gamma_{{\cal O}_{S_f}}]:{\widetilde W}(f)\rightarrow L_{k-1}^p(S^2, \wedge^{0, 1})$ be the corresponding maps under the above trivialization. 
Then 
 $[\gamma_{S_f}](h)=\gamma$ for  $h\in S_f$ and $ [\gamma_{{\cal O}_{S_f}}](h)=T(h)^*(\gamma)$ for $h\in {\widetilde  W}(f).$ 

Hence $[\gamma_{{\cal O}_{S_f}}]$  is smooth , and $[\gamma_{{\cal O}_{S_f}}]$ have the same degree of the smoothness as $T$ has by the following lemma.

\begin{lemma}
For $\gamma \in C^{\infty}(S^2,\wedge^{0, 1} )$, the orbit map  $\Psi_{\gamma}:
G\rightarrow L_{k-1}^p(S^2, \wedge^{0, 1})$ defined by $\Psi_{\gamma}(\phi)=(\phi)^*(\gamma)$ is of class $C^{\infty}.$	
\end{lemma}

The $G$-equivariant extension  $(\beta_i t_{i\nu})_{{\cal O}_{S_f}}$ of $\beta_i t_{i\nu}$ is defined by $(\beta_i t_{i\nu})_{{\cal O}_{S_f}}(h)=:
(\beta_i)_{{\cal O}_{S_f}}(h)\cdot (t_{i\nu})_{{\cal O}_{S_f}}(h)=\beta_i\circ h\cdot t_{i\nu}\circ h$ for $h\in {\widetilde  W}(f)$. It follows from the definition that $(\beta_i t_{i\nu})_{{\cal O}_{S_f}}$ is $G$-equivariant.
We need to show that it is smooth.	

\begin{pro}
Let ${\cal R}={\widetilde  W}(f)\times L_{k-1}^p(S^2, {\bf R})\rightarrow {\widetilde  W}(f)$ be the trivial bundle.  Given a smooth function $\beta:M\rightarrow {\cal R}$, the section $X_{\beta}:{\widetilde  W}(f)\rightarrow {\cal R}$	defined by  $X_{\beta}(h)=\beta\circ h$ is of class $C^{\infty}$. 
\end{pro}

If  we replace ${\bf R}$ by ${\bf C}$,  the same conclusion holds.
This proves the smoothness of $(\beta_i)_{{\cal O}_{S_f}}$. 

Denote  the section $(t_{i\nu})_{{\cal O}_{S_f}}$ by $T_{i\nu}$ for short. Instead of just proving the smoothness of  $T_{i\nu}$, we introduce  a further sheaf theoretic localization of the trivialization  ${\widetilde T}_{k-1,p}|_{ {\widetilde W}(f)}$ by considering the corresponding sheafication,  that will facilitate the proofs of the main theorems  in  this paper and its sequels. For our purpose here, we will only describe  the sections of the sheafication
 on  open sets $D$ containing in some $D_i$ of the above covering   $\{D_i, i=1, \cdot, l\}$.
 To explain this, note   that
 the standard constructions  of  the coordinate charts for ${\widetilde {\cal B}}_{k,p}$  and bundle structure  of  ${\widetilde {\cal T}}_{k-1,p}\rightarrow {\widetilde {\cal B}}_{k,p}$    can obtained 
 as an application of the theory of section functor on the category of FVB (vector bundles with fiber bundle morphisms ) developed by Palais in [P].  The functorial 
 nature of the  section functor  gives rise the desired sheafication. More concretely, for $D\subset D_i$,  let ${\widetilde  W}(f; D)$ be the space of 
 $L_k^p$-maps $h:D\rightarrow V_i\subset M$  with  $||(h-f)|_{D}\|_{k, p}$ less than the prescribed small $\epsilon$ and 
 ${\widetilde  {\cal T}}_{k-1, p}(f; D)\rightarrow {\widetilde  W}(f; D) $  be the  corresponding bundle  defined by  $({\widetilde  {\cal T}}_{k-1, p}(f; D))_h=L_{k-1}^p(D, h^*TM)$ for $h\in {\widetilde  W}_{k, p}(f; D)$.
 For our purpose here, it is sufficient to consider $D=D_i, i=1, \cdots, l.$ Then  we get  the  functorial system  of  the  bundles ${\widetilde  {\cal T}}_{k-1, p}(f; D_i)\rightarrow {\widetilde  W}(f; D_i) $ with respect to the covering.  Each bundle ${\widetilde  {\cal T}}_{k-1, p}(f; D_i)\rightarrow {\widetilde  W}(f; D_i) $, however, is considered as a  morphism   between the two sheaves and 
 the sections $T_{i\nu}, \nu=1, \cdots, m$  before  become   not just   
 the sections of the bundle ${\widetilde  {\cal T}}_{k-1, p}(f; D_i)\rightarrow {\widetilde  W}(f; D_i) $ but also  morphisms between the two sheaves. Moreover, 
 these sections  form a frame of the bundle
 so that they give rise a trivialization of the  bundle. Similarly, the usual process of the local trivialization for the bundle ${\widetilde  {\cal T}}_{k-1, p}$    using the induced parallel transport by the $J$-invariant connection on $TM$  carries over here and produce the "standard''  local trivialization for ${\widetilde  {\cal T}}_{k-1, p}(f; D_i)\rightarrow {\widetilde  W}(f; D_i) $.

{\bf Note:} One of the reasons that we have  brought out the sheaf theoretic aspect in above discussion  is that while sections like  $T_{i\nu}$ or the "constant" sections of  the local bundle ${\widetilde  {\cal T}}_{k-1, p}(f; D_i)\rightarrow {\widetilde  W}(f; D_i) $ are morphism of the corresponding sheaves, the $G$-equivariant extensions  used the first construction in  [L]  mentioned  before can not be interpreted as such  morphisms.

 \begin{theorem}
 The two local trivializations of ${\widetilde  {\cal T}}_{k-1, p}(f; D_i)\rightarrow {\widetilde  W}(f; D_i) $ are $C^{\infty}$ equivalent.  	
 \end{theorem}

The smoothness of the section $(\beta_i t_{i\nu}\circ
f)_{{\cal O}_{S_f}}$ then follows from the above proposition and theorem.
This proves the following main theorem of this paper.

\begin{theorem}
	Given a smooth stable map $f$ and a smooth section 
	
	$\xi\in L_{k-1}^p(S^2, \wedge^{0,1}(f^*TM))$, its "geometric" $G$-equivariant extension $\xi_{{\cal O}_{S_f}}$ above has  the same degree of smoothness as $T$ has. Hence it is either of class $C^m_0$ or of class $C^{\infty}$ accordingly.
\end{theorem}

\medskip
\noindent
{\bf Note:} The  $C^m_0$-smoothness of $\xi_{{\cal O}_{S_f}}$ was already proved in [L] for the local slice $S_f=S_f^E$ obtained by evaluation maps. The proof of this 
is clarified in this paper.

This paper is organized as follows.

In Sec. 2, we give a proof of the  theorem on the smoothness of the $L_k^p$-section 
functor, ${\Gamma}_{k, p} $, applying to a smooth  but possibly non-linear 
 bundle map.  This theorem is part of the theory of Palais in [P] on the Banach space valued section functor on the category of FVB (vector bundles with fiber bundle morphisms). As applications of the theorem, we prove the $C^{\infty}$ smoothness of the coordinate transformations of ${\cal M}ap_{k, p}$ and smoothness of  transition   functions of the bundle ${\widetilde {\cal L} }_{k-1, p}\rightarrow {\cal M}ap_{k, p}$ as well as several related results.
The sheaf theoretic localization of the local trivialization of ${\widetilde {\cal L} }_{k-1, p}|_{{\widetilde W }(f)}\rightarrow {\widetilde W }(f)$ mentioned above  is also discussed in this section.

In Sec. 3, we give the proofs for the two versions of the main theorem.

In Sec. 4, we  consider  an embedding $M\rightarrow {\bf R}^m$ and its induced 
 embedding ${\cal M}ap_{k, p}(S^2, M) =L_k^p(S^2, M)\rightarrow L_k^p(S^2, {\bf R}^m).$ Then  we    construct  a smooth section $\Delta$ of  an infinite dimensional "obstrction" bundle ${\cal N}\rightarrow {\cal T}$.  Here  ${\cal T}$ is an tubular neighborhood of ${\cal M}ap_{k, p}(S^2, M) $ in $L_k^p(S^2, {\bf R}^m)$. The mapping space  ${\cal M}ap_{k, p}(S^2, M) $ then is realized 
 as the zero locus of  the section  $\Delta.$
 Using this construction, we outline a different approach to GW theory. The full details of this part will be treated in a separate paper. 

In Sec. 5, we present  simpler proofs for the smoothness of $L_k^p$-norm and related results.

% \subsection{Applications to GW theory}	

% ${\bullet}$ ${\bullet}$

    \section{ Coordinate Charts and Local Trivializations}
     
     We start with some basic definitions in order to fix our notations.
     
     Let $(M, \omega)$ be a compact  symplectic manifold with a
     symplectic form $\omega$.   Fix   an  $\omega$-compatible almost complex structure $J$.
     Denote  the 
     induced Riemannian metric by $g_J(-, J-).$ Consider the Riemann sphere  $\Sigma=(S^2, i_0)$ with the  standard complex structure $i_0$ and  round metric. Assume that all the geometric data above are of class $C^{\infty}.$
     Let ${\cal  M}ap=:{\cal  M}ap_{k, p}$ be the space of $L_k^p$-maps from $S^2$ to $M$ and  ${\widetilde {\cal  L}}=:{\widetilde {\cal  L}}_{k-1, p}\rightarrow {\cal  M}ap_{k, p}$ be the bundle defined by  ${\widetilde {\cal  L}}_f=L_{k-1}^p(S^2, \wedge^{0,1}_{i_0, J}(f^*(TM)))$ for $f\in {\cal  M}ap_{k, p}.$  We will denote $\wedge^{0,1}_{i_0, J}(f^*(TM))$ by $\wedge^{0,1}(f^*(TM))$ for short.
     
     %Throughout this paper we will always assume that  $p$ is an positive  even integer.
   For the discussion of this section, we  assume further that  $k-2/p>1$.  
     %It is   corresponding  to the following two cases: (i) $k\geq 3$ with $p=2$ or (ii) $k\geq 2$ with $p\geq 4.$
     It follows   that (a) the space $L_{k-1}^p(S^2, {\bf R})$ of $L_k^p$-functions on $S^2$ is a Banach algebra; (b) any element in  ${\cal  M}ap_{k, p}$ is at least of class  $C^1$. Set $m_0=[k-2/p]\geq 1.$

     It is well known that ${\cal  M}ap$  is a smooth Banach manifold and ${\widetilde {\cal  L}}_{k-1, p}\rightarrow {\cal  M}ap_{k, p}$ is a smooth Banach bundle. There are two ways to construct the local trivializations  of the bundle   ${\widetilde {\cal  L}}_{k-1, p}\rightarrow {\cal  M}ap_{k, p}$.
     The first is the standard  one  given in the original paper [G] of Gromov, by identifying the fibers using the induced $J$-invariant parallel transport.
      The other can be obtained  by considering the tangent bundle $T{\cal  M}ap_{k-1, p}$ first, 
     then the simple relation  between   $T{\cal  M}ap_{k-1, p}$ and ${\widetilde {\cal  L}}_{k-1, p}$ gives rise the local trivialization for ${\widetilde {\cal  L}}_{k-1, p}$.  The bundle structure on  $T{\cal  M}ap_{k-1, p}$  is given by Floer in [F1].
     The smoothness of the transition functions between the local trivializations
     of  $T{\cal  M}ap_{k-1, p}$ in [F1]  is a consequence of the smoothness of the coordinate  transformations of  ${\cal  M}ap_{k-1, p}$. On the other hand, a complete
      proof of the smoothness of the transition functions between the local trivializations in [G] as well as the  smoothness of the transition functions between the above two types of local trivializations for  $T{\cal  M}ap_{k-1, p}$  are not presented in the literature of  the GW-theory.   In addition to these local trivializations, in this  paper 
    we  will  introduce a sheaf theoretic type of local trivialization.
    As mentioned before, to prove the main theorem of this paper, it is crucial to show that all these   local trivializations for $T{\cal  M}ap_{k-1, p}$ are $C^{\infty}$-equivalent.
    The main technique tool  to prove such an equivalence as well as other results of similar nature in this paper is the work of Palais in [P] on the category of vector bundles with fiber bundle (hence, possibly  non-linear )morphisms  
    (FVB) and the Banach space valued  section functors on such a category.
    In fact, large  part of the  Floer's work in [F1] can be obtained  as an application of the work in  [P].

     The key ingredient of the construction  in  [P] and [F1] is the following theorem that will be used repeatedly in this paper .
     
     \begin{theorem} ([F1], [P])
     Given a compact smooth manifold $M$, let $B_i\rightarrow M, i=1, 2$ be the bundle of unit ball of the  smooth vector bundle $E_i\rightarrow M$, and $f:B_1\rightarrow B_2$ be a possibly nonlinear smooth bundle map. Assume that $k-\dim (M)/p>0$ and that the space $L_k^p(M)$ of $L_k^p$-functions on $M$ is a Banach algebra.  Then the induced map on $L_k^p$-section, $f_{*}=\Gamma_{k, p} (f):L_k^p(M, B_1)\rightarrow L_k^p(M, B_2)$ is of class $C^{\infty}.$ 
     
     In fact, let $df:B_1\rightarrow L(E_1, E_2)$ be the fiberwise
     derivative of $f$, where $L(E_1, E_2)\rightarrow M$ is the vector bundle with the fiber $L(E_1, E_2)_m=L((E_1)_m, (E_2)_m)$ for $m\in M.$ Then $df$ is a smooth  bundle map and $$D(f_{*})=D(\Gamma_{k, p} (f))=\Gamma_{k, p} (df):
     \Gamma_{k, p}(B_1)=L_k^p(M, B_1)\rightarrow  \Gamma_{k, p}(L(E_1, E_2))=$$ $$L_k^p(M, L(E_1, E_2))\subset L(\Gamma_{k,p}( E_1), \Gamma_{k,p}( E_2)
    )=L(L_k^p(M, E_1), L_k^p(M, E_2)).$$

     \end{theorem}
     
     \proof 
     
     We only  outline a proof. For more details, see [F1] and [P].
     Clearly
      we only need to show that (1) $f_{*}=\Gamma_{k, p} (f)$ is continuous and (2) $D(f_{*})=D(\Gamma_{k, p} (f))$ is equal to $\Gamma_{k, p} (df)$.
     Indeed, since the derivative $df$ along the fiber  is a smooth bundle map so that it plays the same role as $f$ does, (1) and (2) imply that $f_{*}$ is of class $C^{\infty}$ by induction.

     ${\bullet}$ Proof of (1):
     
     Let $s_1, s_2\in L_k^p(M, B_1)$. Then for $j\leq k$, $$\|D ^j(f_{*}(s_1)-f_{*}(s_2))\|_{0, p}=\|D ^j(f\circ s_1-f\circ s_2)\|_{0, p}$$ $$=\int_{0}^{1} \|D ^j(\frac {df}{dt}(s_2+t(s_1-s_2)))\|_{0, p} dt$$ 
     
     is bounded by terms $\|df\|_{C^k}\cdot \|D ^i(s_1-s_2)\|_{0, p}$ with $i\leq j.$ In fact, a better bound is $\|df\|_{k, q}\cdot \|D ^i(s_1-s_2)\|_{k, p}$.
     Now 
     for $p\geq 2$,  $q\leq 2$, $L_k^p( S^2)$ is continuously embedded into $L_k^q( S^2)$ so that we have the bound  $\|df\|_{k, p}\cdot \|D ^i(s_1-s_2)\|_{k, p}$.

      ${\bullet}$ Proof of (2):
     
     For $s\in  L_k^p(M, B_1)$,  $\xi\in  L_k^p(M, E_1) $ and $j\leq k,$
     $$\|D ^j(f_{*}(s+t\xi)-f_{*}(s)-\Gamma_{k, p}(df)_s(t\xi))\|_{0, p}=\|D ^j(f\circ (s+t\xi)-f\circ s- tdf_s\circ \xi)\|_{0, p}$$ $$ = \|D ^j(\int_0^1 (df (s+\nu t\xi)-df(s))d\nu (t\xi))\|_{0, p}\leq \int_0^1 \|D ^j( (df (s+\nu t\xi)-df(s) )(t\xi))\|_{0, p}d\nu $$  $$\leq \int_0^1 \int_0^1\|D ^j( (d^2f (s+\mu\nu t\xi)(\nu t\xi) )(t\xi))\|_{0, p}d\nu d \mu ,$$
     which is bounded by 
      $ \|df^2\|_{C^k}\cdot \|t\xi\|^2_{k, p}$ under the assumption that
      $L_k^p(M)$ is a Banach algebra. As remarked above, a better bound is 
     $ \|df^2\|_{k, p}\cdot \|t\xi\|^2_{k, p}$
     
     \QED

    {\bf Note:} Using the better bounds in  the  proof above, we  only to  assume that (a )the bundle $E_i\rightarrow M, i=1,2$ and the bundle  map  $f:B_1\rightarrow B_2$ are of class $L_k^p$ (or $C^k$); (b)   $f$ is of class $C^{\infty}$ along the fiber such that the derivatives $df$ and $d^2f$ along the fiber as the corresponding  bundle  maps are of class  $L_k^p$ (or $C^k$). Then the above theorem is still  true.

    We now apply above theorem to construct the standard coordinate transformations
    for  ${\cal M}ap$ and the transition functions for ${\widetilde {\cal L}}.$
   In the following unless specified otherwise, we will  assume that the center $f$ of each coordinate is of class $C^{\infty}.$ However, one can verify that the conditions  in  the above note  are satisfied for the discussion below so that the same results are true when the center $f\in  {\cal M}ap_{k, p}$ is  of class $L_k^p.$
   
   \noindent
   ${\bullet }$   ${\bullet }$ Smoothness of the coordinate transformations  of  ${ {\cal M}}ap_{k, p}$.

   Recall the definition of  natural coordinate chart ${\widetilde W}(f)$ for $f\in { {\cal M}}ap_{k, p}$.
     Consider the  smooth bundle
     $E=f^*(TM)\rightarrow S^2$. Denote its     sub-bundle of unit balls by $B\rightarrow \Sigma=S^2.$   
   Let ${\hat W}= {\hat W}_{\epsilon}$ be the
   $\epsilon$-ball of $L_k^p(S^2, B)\subset L_k^p(S^2, f^*(TM)).$ Then the  $\epsilon$-neighborhood of $f$ in ${\cal M}ap$ is defined to be   ${\widetilde W}(f)=:Exp_f ({\hat W})$, where $(Exp_f\xi)(x)=exp_{f(x)}\xi(x)$. 
   %Here the exponential map is taken with respect to an 
   %$f$-dependent metric that will be specified  later.  
   
    Given two smooth elements $f_i, i=1,2 \in {\cal M}ap$, the above theorem implies that  the coordinated transformation $\Psi_{21}:{\widetilde W}_{f_1}\rightarrow {\widetilde W}_{f_2}$
    is of class $C^{\infty}.$ Indeed,  assume $f_1$ and $f_2$ are $C^{0}$-close each other so that the  map  $\psi_{21}:B_1\rightarrow B_2$ given by $\psi_{21}(b) = exp^{-1}_{f_2(x)}\circ exp_{f_1(x)}(b)$ is well-defined on $B_1$, where $b\in (B_1)_x$. Then   $\psi_{21} $ is of class $C^{\infty},$ and the 
    coordinated transformation $\Psi_{21}:{\widetilde W}_{f_1}\rightarrow {\widetilde W}_{f_2}$ is  the restriction of $(\psi_{21})_{*}=\Gamma_{k, p}(\psi_{21})$ to ${\widetilde W}({f_1})$. Hence by the above theorem  $\Psi_{21}$ is  of class of $C^{\infty}$.
     
     \noindent
     ${\bullet }$   ${\bullet }$ Smoothness of the  transition functions  of  ${\widetilde {\cal L}}.$
     
    Next recall  the local trivialization of the bundle  ${\widetilde {\cal L}}={\widetilde {\cal L}}_{k-1, p}$ over the local chart  ${\widetilde W}_f$.
    For any $h\in  {\widetilde W}_f$ and $\xi\in {\widetilde {\cal L}}_f$,  the local trivialization ${\Pi}_f:{\widetilde W}_f\times {\widetilde {\cal L}}_f\rightarrow 
    {\widetilde {\cal L}}|_{{\widetilde W}_f}$ is given by ${\Pi}_f(h, \xi)(x)=P_{h(x)f(x)}\xi(x)$ for $x\in S^2$  where  $ P_{h(x)f(x)}$ is the induced action  by  the $J$-invariant parallel transport along the shortest geodesic from $f(x)$ to $h(x).$
    
    Given $f_1$ and $f_2$ in  ${\cal M}ap$
     as above,  the  transition function between  two local trivializations then is defined by ${\Pi}_{21}=:{\Pi}^{-1}_{f_2}\circ {\Pi}_{f_1}:{\widetilde W}_{f_1}\times {\widetilde {\cal L}}_{f_1}\rightarrow 
     {\widetilde W}_{f_2}\times {\widetilde {\cal L}}_{f_2}$. To see that ${\Pi}_{21}$ is of class $C^{\infty}$, we make some reductions.
     
     Note that any $\xi\in  {\widetilde {\cal L}}_f$  is a summation of the elements  $\eta \otimes \gamma$ with $\eta\in L_{k-1}^p(S^2, f^*TM)$ and $\gamma\in C^{\infty}(S^2, \wedge^{0, 1} ).$ In other words,  ${\widetilde {\cal L}}_f=L_{k-1}^p(S^2, f^*TM)=C^{\infty}(S^2, \wedge^{0, 1} )\otimes L_{k-1}^p(S^2, \wedge^{0, 1}(f^*TM))$. Note that identification of the fibers induced by the parallel transport  in local trivialization for ${\widetilde {\cal L}}$ only acts on $L_{k-1}^p(S^2, f^*TM).$
     This implies that  ${\widetilde {\cal L}}\simeq {\bf \Omega}^{0, 1}\otimes {\widetilde {\cal T}}_{k-1, p}$. Here ${\bf \Omega}^{0, 1}\rightarrow {\widetilde {\cal B}}$ is the  trivial bundle  defined by ${\bf  \Omega}^{0, 1}_f=C^{\infty}(S^2, \wedge^{0, 1})$, and  the bundle ${\widetilde {\cal T}}_{k-1, p}$ with the fiber  ${\widetilde {\cal T}}_{k-1, p}|_{f}=L_{k-1}^p(S^2, f^*TM)$  has the local trivialization 
     induced by the parallel transport as for ${\widetilde {\cal L}}.$
     Thus we only need to show that the corresponding transition function ${\Pi}_{21}$ for ${\widetilde {\cal T}}_{k-1, p}$ is of class $C^{\infty}.$

    % We now  prove that the  transition function ${\Pi}_{21}$ of ${\widetilde {\cal T}}_{k-1, p}$  is of class $C^{\infty}.$

     To this end,  let  $p_{i}:B_i\subset E_i\rightarrow S^2$ 
     be the  unit ball bundle  of $E_i=:f_i^*TM$ before. Consider   the pull-back bundles 
     ${\bf E}_i= (Exp_{f_i})^*(TM)\rightarrow B_i, i=1,2 $  by the maps
     $Exp_{f_i}:B_i\subset E_i\rightarrow M.$ 
      Then the fiber $({\bf E}_i)_{b}$ at $b\in B_i$ with $p_{i}(b)=x$ is $T_{exp_{f(x)}b}M.$
     
     Next consider the pull-back  of the bundles $E_i\rightarrow S^2$ by $p_i:B_i\rightarrow S^2$,  $p_i^*(E_i)\rightarrow B_i.$  The  bundle map $\pi_{f_i}: p_i^*(E_i)
      \rightarrow {\bf E}_i$ is defined by $\pi_{f_i}(b, e)=(b, P({exp_{f_i(x)}b, f_i(x)})( e))$  for $(b, e)\in p_i^*(E_i)$ with $b\in B_i, p_i(b)=x$ and $ e\in (E_i)_x$, where $P({exp_{f_i(x)}b, f_i(x)})( e)$ is the parallel transport of $e$ along the geodesic from $f_i(x)$  to $exp_{f_i(x)}b$ by  a fixed $J$-invariant connection on $(TM, J)$. %Note that  here the geodesic is  with respect to the  $f_i$-dependent metric. 
      It is easy to see that $\pi_{f_i}$  is of class $C^{\infty}$.
     Then  the bundle map $\pi_{21}=:\pi_{f_2}^{-1}\circ \pi_{f_1}:p_1^*(E_1) \rightarrow p_2^*(E_2)$ is of class $C^{\infty}$ as well.

     In order to applying the above theorem, we need to find the corresponding maps  between  the relevant  vector bundles on $S^2$ rather than on $B_i$ as above.
     To this end,  consider the bundle $H=H_{21}=Hom(E_1, E_2)\rightarrow S^2$ and the bundle map
     $T:B_1\rightarrow H$ given by $T(b)(e)=\pi_{21}(b, (p_1^*(e))(b))$. Here $p_1^*(e)$ is  the  "vertical" lifting of $e\in (E_1)_x$, which is a  section of 
     $p_1^*(E_1)$ along the fiber $(B_1)_x$ and $(p_1^*(e))(b)$ is  its value at $b\in (B_1)_x.$ Then the smoothness of $\pi_{21}$ implies that $T$ is of class
     $C^{\infty}.$
     
     Now applying the $L_{k-1}^p$-section  functor ${\Gamma}_{k-1, p}$, we get
     Banach spaces $ {\Gamma}_{k-1, p}(H)$ and $ {\Gamma}_{k-1, p}(E_i), i=1,2.$ 
     The paring $<-, ->:{\Gamma}_{k-1, p}(H)\times 
     {\Gamma}_{k-1, p}(E_1)\rightarrow {\Gamma}_{k-1, p}(E_2)$ is of class $C^{\infty}$.  By the above theorem, the smoothness of $T$  implies that 
     ${\Gamma}_{k-1, p}(T):{\Gamma}_{k, p}(B_1)\rightarrow {\Gamma}_{k, p}(H)\rightarrow {\Gamma}_{k-1, p}(H)$ is of class $C^{\infty}.$
     
     Consequently the map ${\bf T}=:<-, ->\circ  ({\Gamma}_{k-1, p}(T)\times Id): {\Gamma}_{k, p}(B_1)\times {\Gamma}_{k-1, p}(E_1)\rightarrow  {\Gamma}_{k-1, p}(H)\times {\Gamma}_{k-1, p}(E_1)\rightarrow 
    {\Gamma}_{k-1, p}(E_2)$ is of class $C^{\infty}.$
    
    Now $\Pi_{21}(\xi, \eta)=(\Psi_{21}(\xi), {\bf T}(\xi, \eta ))$. This implies that
    $\Pi_{21}$ is of class $C^{\infty}.$
    
     %Now  for each $i$,   we have  a  bundle sequence of length two    ${ S}_i=:p_i^*(E_i^{0,1})\rightarrow B_i\rightarrow S^2, i=1, 2, $. The two bundle maps $(\pi_{21}, \psi_{21})$ above together define an isomorphism between these two 
     %bundle sequences.  

  % Then    $\Gamma_{k, k-1;p}(S_i):=\Gamma_{k, p}( B_i)\times \Gamma_{k-1, p}(p_i^*(E_i^{0,1})) $ and $\Gamma_{k, k-1;p}(\psi_{21}, \pi_{21})=\Gamma_{k, p}(\psi_{21})\times \Gamma_{k-1, p}(\pi_{21}).$
    % One can show that in this situation, the theorem above is still true so that $\Gamma_{k, k-1;p} (\psi_{21}, \pi_{21})$  is of class $C^{\infty}$.

      %Now it follows from the construction that ${\Pi}_{21}=\Gamma_{k, k-1;p} (\psi_{21},\pi_{21})$ so that  it is of class  $C^{\infty}$.

    %It contains an  open submanifold ${\widetilde {\cal B}}$ of stable $L_k^p$-maps defined as follows.
     
     ${\bf \bullet}{\bf \bullet}$ Smooth equivalence between the local trivializations of $T{\widetilde {\cal B}}_{k-1, p}$ and ${\widetilde {\cal T}}_{k-1, p}$

      Note that if we replace $k-1$ by $k$, the bundle ${\widetilde {\cal T}}_{k, p}$ is just the ($C^{\infty}$) tangent bundle $T{\widetilde {\cal B}}_{k, p}$ but with different local trivializations. 
      The local trivializations for  $T{\widetilde {\cal B}}_{k, p}$ induce  the local trivializations  and hence 
      a $C^{\infty}$ bundle structure for  $T{\widetilde {\cal B}}_{k-1, p}$.
      We now show that the set theoretic identification $T{\widetilde {\cal B}}_{k-1, p}\rightarrow {\widetilde {\cal T}}_{k-1, p}$ is in fact a bundle isomorphism
       of class $C^{\infty}.$ We will prove this by showing that the transition functions between of the local trivializations of $T{\widetilde {\cal B}}_{k-1, p}$ and ${\widetilde {\cal T}}_{k-1, p}$ are of class $C^{\infty}.$
       
       First recall the local trivialization of $T{\widetilde {\cal B}}_{k-1, p}\rightarrow {\widetilde {\cal B}}_{k, p}$ over a local  chart
        ${\widetilde {W}}(f)$, $D_f=:D({Exp}_f):{\widetilde W}(f)\times (T{\widetilde {\cal B}})_f\rightarrow 
        T{\widetilde {\cal B}}|_{{\widetilde W}_f}$  is
         given by ${D}_f(h, \xi)(x)=(D({Exp}_{f(x)}))_{h(x)}(\xi(x))$ for $x\in S^2.$   To show that the transition function between $D_f$ and $\Pi_f$  is of class $C^{\infty}$, we need to find the corresponding bundle maps as  we did for the transition function $\Pi_{21}.$
        
        Consider the bundle $p:E=f^*(TM)\rightarrow S^2$,  its unit ball bundle
        $p:B\rightarrow S^2$,  the pull-back bundles $p^*(E)\rightarrow B$ by $p$
        and ${\bf E}=: (Exp_{f})^*(TM)\rightarrow B  $  by 
        $Exp_{f}:B\subset E\rightarrow M.$ 
        Define the bundle isomorphism $I=I_f:p^*(E)\rightarrow {\bf E}$ by $I(b, e)=(b,(D({Exp}_{f(x)}))_{{Exp}_{f(x)}^{-1}(b)}(e ))$ where $x=p(b).$  Then $I$ is of class $C^{\infty}.$ Let $A=A_f=:\pi_f^{-1}\circ I_f:p^*(E)\rightarrow p^*{ (E)}$ be the bundle automorphism. Then $A$ is of class $C^{\infty}$.
        As before  define  a $C^{\infty}$ bundle  map $T:B\rightarrow H$  by $T(b)(e)=A(b, (p^*(e))(b)),$ where 
       the bundle $H=Hom(E, E)\rightarrow S^2$.  Here $p^*(e)$ is  the section the fiber $(B_1)_x$ of the  vertical lifting of $e\in (E)_x$ and $(p^*(e))(b)$ is the its value at $b\in (B_1)_x.$  Then the argument before implies that
       ${\Gamma}_{k-1, p}(T):{\Gamma}_{k, p}(B)\rightarrow {\Gamma}_{k, p}(H)\rightarrow {\Gamma}_{k-1, p}(H)$ is of class $C^{\infty}$ so that 
       the map ${\bf T}=:<-, ->\circ  ({\Gamma}_{k, p}(T)\times Id): {\Gamma}_{k, p}(B)\times {\Gamma}_{k-1, p}(E)\rightarrow  {\Gamma}_{k-1, p}(H)\times {\Gamma}_{k-1, p}(E)\rightarrow 
       {\Gamma}_{k-1, p}(E)$ is of class $C^{\infty}.$

       Now the transition function ${\bf A}$ between the two local trivializations
       $D_f$ and $\Pi_f$ above is given by   ${\bf A}(\xi, \eta)=(\xi, {\bf T}(\xi, \eta ))$. This implies that
       ${\bf A}$ is of class $C^{\infty}$, hence proves the following proposition.
       
       \begin{pro}
       	 The transition functions between the local trivializations of ${\widetilde {\cal T}}_{k-1, p}\simeq T{\widetilde {\cal B}}_{k-1, p}$ are of class  $C^{\infty}$.
       \end{pro}

      \begin{cor}
      	The section $s=:{\bar{\partial}}_J:{\widetilde {\cal B}}_{k, p}\rightarrow :{\widetilde {\cal L}}_{k-1, p}$  defined by $ s(f)=df+J(f)\cdot df\circ i_0$ is 
      	of class $C^{\infty}.$
      	
      \end{cor} 
      
      \proof 
       
       The result is  proved in [F]  and [P]  for $s:{\widetilde {\cal B}}_{k, p}\rightarrow :{\widetilde {\cal L}}_{k-1, p}$ with respect to the smooth structure on ${\widetilde {\cal L}}_{k-1, p}$ induced from $T{\widetilde {\cal B}}_{k-1, p}$. 
       
       \QED
        
       ${\bf \bullet} $ Further sheaf-theoretic localization of a trivialization of the bundle ${\widetilde {\cal T}}_{k-1, p}$

      \medskip
       To facilitate the proof of the main theorem,  we now introduce a   sheaf-theoretic localization of certain sections of  the   local  bundle ${\widetilde {\cal T}}_{k-1, p}|_{{\widetilde { W}}(f)}\rightarrow {\widetilde {\cal W}}(f)$.

        Let $\cup_{i=1}^lD_i=S^2$ be an fixed open covering of $S^2$ and $\alpha_i, i=1, \cdots, l$ is a partition of unit subordinated to the covering.
        Each $\alpha_i$ induces a bundle morphism $E_{\alpha_i}:{\widetilde {\cal T}}_{k-1, p}\rightarrow {\widetilde {\cal T}}_{k-1, p}$ defined by $E_{\alpha_i}(\eta)=\alpha_i\cdot \eta$ for $\eta \in {\widetilde {\cal T}}_{k-1, p}|_h=L_{k-1}^p(S^2, h^*TM)$  with  $h\in {\widetilde {\cal B}}_{k, p}$. 
      It is easy to see that $E_{\alpha_i}$ is 
        smooth. Indeed in the any of above local trivialization ${\widetilde {\cal T}}_{k-1, p}|_{{\widetilde { W}}(f)}\simeq {{\widetilde { W}}(f)}\times L_{k-1}^p(S^2, f^*TM)$, the corresponding map $[E_{\alpha_i}]_f:
        {{\widetilde { W}}(f)}\times L_{k-1}^p(S^2, f^*TM)\rightarrow {{\widetilde { W}}(f)}\times L_{k-1}^p(S^2, f^*TM)$  given by 
        $[E_{\alpha_i}]_f(h, \xi)=(h, \alpha_i\cdot \xi)$ is "linear" and continuous, hence smooth.
        
        Thus given a section $s:{\widetilde {\cal B}}_{k, p}\rightarrow {\widetilde {\cal T}}_{k-1, p},$ $s$ is of class $C^{k}, k=0, 1, \cdots, \infty,$ if and only  if $E_{\alpha_i}(s)$ is  for all $ i=1, \cdot, l.$
        
        Now let $\cup_{i=1}^l V_i=M$ be a fixed open covering.
        % Denote the covering $\{V_i, i=1, \cdot, l\} $ by ${\it V}.$
         Assume that the cover
        ${\it D}$ of $S^2$  above  satisfies the condition that   $h(D_i) \subset V_i, i=1, \cdots, l$ for any $h\in  {{\widetilde { W}}(f)}$. In particular, 
        for $(h, \phi)\in S_f\times G_e$, $h\circ \phi (D_i) \subset V_i.$
        Let $V_i\subset \subset V_i', i=1, \cdots, l$ and 
        $\beta_i, i=1, \cdots, l$ be  cut-off functions defined on $V'_i$   such that 
        $\beta_i=1$ on $V_i$.  For each $i$, let $\alpha'$ be a cut-off function on $S^2$ supported on $D'_i$  with $D_i\subset \subset D_i'$ and $\alpha'_i=1$ on $D_i$.
        Fix a smooth local  frame ${\bf t}_{i}=\{t_{i1}, \cdots, t_{im}\}$ of $TM$ on
         $V_i'.$ Then 
         for $h\in {{\widetilde { W}}(f)}$ and $x\in D_i$, 
        $$(E_{\alpha_i}(s))(h)(x)=\Sigma_{\nu}(s_i^{\nu}(h))(x)\cdot t_{i\nu}(h(x))  =\Sigma_{\nu}(s_i^{\nu}(h))(x)\cdot \beta_i (h(x)) t_{i\nu}(h(x))
        $$ $$ =\Sigma_{\nu}(s_i^{\nu}(h))(x)\cdot \alpha'_i (x) t_{i\nu}(h(x))
        .$$
       
       Note that by our assumption for any $h\in {{\widetilde { W}}(f)},$
       (1 )the pull back frame $h^*({\bf t}_{i})$ is defined over $h^{-1}(V_i)$ that contains $D_i$ hence the support of $(E_{\alpha_i}(s))(h)$; (2)the pull back  of the global sections sections $\beta_i{\bf t}_{i}=:\{\beta_i{ t}_{i1}, \cdots, \beta_i{ t}_{im}\}$, $h^*(\beta_i{\bf t}_{i})$    are global sections
        of the bundle $h^*TM\rightarrow S^2$, which plays the role of a "global frame" of the bundle  $h^*TM\rightarrow S^2$ to express the  sections like $(E_{\alpha_i}(s))(h)$  supported on $D_i$; (3) by multiplying by $\alpha'_i$ the frame $h^*({\bf t}_{i})$ above becomes $E_{\alpha'_i}(h^*({\bf t}_{i}))=:\alpha'_ih^*({\bf t}_{i})$ that is also a "global frame" of the bundle  $h^*TM\rightarrow S^2$.
        
        Hence for each $i$, we get three  corresponding collections of sections, 
        denoted by ${\bf T}_i=\{T_{i1}, \cdots, T_{im}\}$, $E_{\beta_i} ({\bf T}_i)=\{E_{\beta_i} (T_{i1}), \cdots, E_{\beta_i} (T_{im})\}$ and $E_{\alpha'_i}({\bf T}_i)= \{E_{\alpha_i} (T_{i1}), \cdots, E_{\alpha_i} (T_{im})\}$. The section $T_{i\nu}: {{\widetilde { W}}(f; D_i)}\rightarrow {\widetilde  {\cal T}}_{k-1, p}(f; D_i)$ is defined by $T_{i\nu}(h)=t_{ij}\circ h$ for $h\in {\widetilde  W}(f; D_i).$ Here 
        ${\widetilde  W}(f; D_i)$ consists  of $L_k^p$-maps $h:D_i\rightarrow V_i\subset M$ such that $||(h-f)|_{D_i}\|_{k, p}$ less than the prescribed small $\epsilon$;  and  
        the bundle ${\widetilde  {\cal T}}_{k-1, p}(f; D_i)$  is defined by the  same  formula as before, $({\widetilde  {\cal T}}_{k-1, p}(f; D_i))_h=L_{k-1}^p(D_i, h^*TM)$ for $h\in {\widetilde  {\cal T}}_{k-1, p}(f; D_i)$. The section $E_{\beta_i} (T_{i\nu}): {\widetilde  W}(f )\rightarrow {\widetilde  {\cal T}}_{k-1, p}|_{{\widetilde W}(f )}$ is defined by $E_{\beta_i} (T_{i\nu})(h)=(\beta_i\cdot t_{ij})\circ h$ for $h\in {\widetilde  W}.$ The section $E_{\alpha_i} (T_{i\nu}): {\widetilde  W}(f )\rightarrow {\widetilde  {\cal T}}_{k-1, p}|_{{\widetilde W}(f )}$ is defined similarly.

        The main theorem is the following.	
        	
        \begin{theorem}
        The section $E_{\alpha_i} (T_{i\nu})\& E_{\beta_i} (T_{i\nu}) : {\widetilde  W}(f )\rightarrow {\widetilde  {\cal T}}_{k-1, p}|_{{\widetilde W}(f )}$ 	
        	as well as the section $ T_{i\nu}: {{\widetilde { W}}(f; D_i)}\rightarrow {\widetilde  {\cal T}}_{k-1, p}(f; D_i)$ are smooth of class $C^{\infty}$.

        \end{theorem}	
        
        To simplify our notations, we will drop the subscript $i$ for the discussion below.
        The theorem will be derived as a corollary of  the smooth equivalence of the following two local trivializations of the bundle ${\widetilde  {\cal T}}_{k-1, p}(f; D)\rightarrow  {{\widetilde { W}}(f; D)}.$ The first local trivialization is just the standard one induced by the $J$-invariant  parallel transport.
        It can be reformulate as follows. Let $ Q_{\nu}:{{\widetilde { W}}(f; D)}\rightarrow {\widetilde  {\cal T}}_{k-1, p}(f; D)$ be the corresponding section obtained as the "constant" extension  of the point-section $t_{\nu}\circ f|_{D}\in ({\widetilde  {\cal T}}_{k-1, p}(f; D))_f=L_{k-1}^p(D, f|_{D}^*{TM})$  by the standard local trivialization induced by the parallel transport along short geodesics. 
        Then the local trivialization $\Pi^1:{{\widetilde { W}}(f; D)}\times L_{k-1}^p(S^2, f_D^*TM)\rightarrow {\widetilde  {\cal T}}_{k-1, p}(f; D)$ is given by $\Pi^1(h, \xi)=\Sigma_{\nu}\xi^{\nu}Q_{\nu}(h)$ for $(h, \xi)\in 
        {{\widetilde { W}}(f; D)}\times L_{k-1}^p(S^2, f_D^*TM)$ with $\xi=\Sigma_{\nu}\xi^{\nu}Q_{\nu}(f|_D),$ where $\xi^{\nu}\in L_{k-1}^p(D, {\bf C})$. Denote the collection of $ Q_{\nu}, \nu=1, \cdots, m$ by $ {\bf Q} (= {\bf Q}_i, i=1, \cdots, l)$. Thus like  ${\bf T}$,  $ {\bf Q}$  is also a frame
        for ${\widetilde  {\cal T}}_{k-1, p}(f; D)$.
        % in the sense that ${\widetilde  {\cal T}}_{k-1, p}(f; D)$ is a free module over the Banach algebra
        %$L_{k-1}^p(D, {\bf R})$ of rank $m=\dim (M).$  
       % Therefore, like the finite dimensional case,  two  point views here,   the trivialization and free module over the  Banach algebra  are equivalent.
        
        The second trivialization $\Pi^2:{{\widetilde { W}}(f; D)}\times L_{k-1}^p(S^2, f_D^*TM)\rightarrow {\widetilde  {\cal T}}_{k-1, p}(f; D)$ is given  similarly by replacing $Q_{\nu}$ above by $T_{\nu}$. Hence  $\Pi^2(h, \xi)=\Sigma_{\nu}\xi^{\nu}T_{\nu}(h)$ for $(h, \xi)\in 
        {{\widetilde { W}}(f; D)}\times L_{k-1}^p(S^2, f_D^*TM)$ with $\xi=\Sigma_{\nu}\xi^{nu}T_{\nu}(f|_D).$
        
        \begin{pro}
        
         The two local trivializations $\Pi^1\& \Pi^2:{{\widetilde { W}}(f; D)}\times L_{k-1}^p(S^2, f_D^*TM)\rightarrow {\widetilde  {\cal T}}_{k-1, p}(f; D)$  are smoothly equivalent.
        
        \end{pro}

        \proof
      
       As  before,  the proposition  will be proved by applying the theory of 
        TVB in [P] by  finding  the corresponding bundle maps.
       In steady of essentially repeating what we did before, we give an "universal"
         treatment that works for the case here as well as the cases before  with necessary   modifications. 
         
       Let $V$ be one of those $V_i, i=1, \cdots, l$ above. Denote  the bundle $TV$ by  $p:E=:TV\rightarrow S^2$,  its unit ball bundle by
       $p:B\rightarrow S^2$. Consider   the pull-back bundles $p^*(E)\rightarrow B$ and ${\bf E}=: (Exp)^*(E)\rightarrow B  $  by the maps $p:B\rightarrow S^2$
       and  
       $Exp=Exp_V:B\subset E\rightarrow M$  respectively.

       Define the following two bundle isomorphisms $I_1\&I_2: p^*(E)\rightarrow {\bf E}$ as follows. For any $(b, e) \in  p^*(E) $  with $b\in B, p(b)=x$ and $ e\in E_x$, $I_1=:p^*(E)\rightarrow {\bf E}$ by $I_1(b, e)=(b, P({exp_xb, x})( e))$, where  $P({exp_xb, x})( e)$ is the parallel transport of $e$ along the geodesic from $x$  to $exp_x b$ by   the  fixed $J$-invariant connection on $(TM, J)$.  Then $I_1$ is of class $C^{\infty}.$ Thus this is the bundle
       isomorphism relevant to the standard trivialization induced from the parallel transport.  The second bundle isomorphisms $I_2: p^*(E)\rightarrow {\bf E}$ is defined by $I_1(b, e)=(b, \Sigma_{\nu}a^{\nu}t_{\nu}(exp_xb))$ for  $b\in B, p(b)=x$ and $ e=\Sigma_{\nu}a^{\nu}t_{\nu}(x)\in E_x.$ Again $I_2$ is of class $C^{\infty}.$

       Let $A=:(I_2)^{-1}\circ I_1:p^*(E)\rightarrow p^*{ E}$ be the bundle automorphism. Then $A$ is of class $C^{\infty}$.
       
       Now let $f_D=f|_D:D\rightarrow V\subset M$. Consider the pull-back bundle $p_f=p_{f_D}: E_f=f_D^*(E)\rightarrow D$, its unit ball bundle   $p_f:
       B_f=f_D^*(B)\rightarrow D$ as well as the pull-backs  of all the related commutative diagrams.  In particular, consider  the pull-back bundle $p_f^*(E_f)\rightarrow B_f$ of $E_f$ by the map $p_f:
       B_f\rightarrow D.$ Then the bundle automorphism $A$ induces a corresponding $C^{\infty}$-bundle automorphism $A_f:p_f^*(E_f)\rightarrow p_f^*(E_f).$

       As before  define  a $C^{\infty}$ bundle  map $T_f:B_f\rightarrow H_f$  by $T(b)(e)=A(b, (p^*(e))(b)),$ where 
       the bundle $H_f=Hom(E_f, E_f)\rightarrow D$.  Here $p^*(e)$ is  the  vertical lifting of $e\in (E_f)_x$  along the fiber $(B_f)_x$   and $(p_f^*(e))(b)$ is the its value at $b\in (B_f)_x.$  Then the argument before implies that
       ${\Gamma}_{k-1, p}(T_f):{\Gamma}_{k, p}(B_f)\rightarrow {\Gamma}_{k, p}(H_f)\rightarrow {\Gamma}_{k-1, p}(H_f)$ is of class $C^{\infty}$ so that 
       the map ${\bf T_f}=:<-, ->\circ  ({\Gamma}_{k, p}(T_f))\times Id: {\Gamma}_{k, p}(B_f)\times {\Gamma}_{k-1, p}(E_f)\rightarrow  {\Gamma}_{k-1, p}(H_f)\times {\Gamma}_{k-1, p}(E_f)\rightarrow 
       {\Gamma}_{k-1, p}(E_f)$ is of class $C^{\infty}.$

       Now the transition function ${\bf A}_f={\bf A}_{f_D}$ between the two local trivializations
       $\Pi^1_f$ and $\Pi^2_f$ above is given by   ${\bf A}_f(\xi, \eta)=(\xi, {\bf T}(\xi, \eta ))$. This implies that
       ${\bf A}_f$ is of class $C^{\infty}$,  hence proves the  proposition.
        \QED

         ${\bf \bullet }{\bf \bullet }$  Proof of the theorem:
      
       Clearly    the smoothness of ${ T_{\nu}}$ follows from the above proposition.
       
       The section $E_{\alpha} (T_{\nu}): {\widetilde  W}(f )\rightarrow {\widetilde  {\cal T}}_{k-1, p}|_{{\widetilde W}(f )}$ considered as a
        map $ {\widetilde  W}(f )\rightarrow L_{k-1}^p(S^2, f^*TM)={\widetilde  {\cal T}}_{k-1, p}|_{f }$ is  a composition $E_{\alpha} (T_{\nu})=E_{\alpha}\circ
        T_{\nu} \circ R: {\widetilde  W}(f )\rightarrow {{\widetilde { W}}(f; D)}\rightarrow  L_{k-1}^p(S^2, f_D^*TM )\rightarrow  L_{k-1}^p(S^2, f^*TM)$. Here $R: {\widetilde  W}(f )\rightarrow {{\widetilde { W}}(f; D)}$ is the restriction map and $E_{\alpha}: L_{k-1}^p(S^2, f_D^*TM )\rightarrow  L_{k-1}^p(S^2, f^*TM)$ is the multiplication by $\alpha$ with support of $\alpha$ containing in $D$. Clearly both $R$ and $E_{\alpha}$ are continuous and linear (in the local chart for $R$), hence smooth so that $E_{\alpha} (T_{\nu})$ is a smooth section.
       
       To see the smoothness of $ E_{\beta} (T_{\nu}) : {\widetilde  W}(f )\rightarrow {\widetilde  {\cal T}}_{k-1, p}|_{{\widetilde W}(f )}$,    recall that  support of ${\beta}$ is contained in $V'$ such that $\beta =1$ on $V\subset\subset V'$. We will assume that the local frame ${\bf t}$ is defined 
       on slightly larger open neighborhood $V^{(4)} $ of ${\bar V''}$ with $V'\subset\subset V'''\subset\subset V''$. Then we may assume that $ f^{-1}({\bar V''})\not =S^2$, otherwise  the image of  any $h\in {\widetilde W}(f )$ is already contained in $V^{(4)}$
       so that $T_{\nu}$  is already a smooth section on ${\widetilde W}(f )$ by above proposition.
       
        Now defined the inverse images  $D'=f^{-1}(V')\subset\subset D'''=f^{-1}(V''')\subset\subset D''=f^{-1}(V'')$.  Note that  by our assumption above, the relations such as ${\bar D'}\subset f^{-1}({\bar V'})\subset D'''$  hold.
        We may assume that
       $D'$ here is the same as the one defined before. Since  
       ${\bar D'} \subset f^{-1}({\bar V'}) \subset  f^{-1} (V''')$, 
       we may assume that  for any $h\in {\widetilde  W}(f )$, the image of $h( f^{-1}({\bar V'}))\subset V'''.$  Indeed since $f( f^{-1}({\bar V'})$ is compact and contained in $V''', $ there exists a small $ \epsilon $-neighborhood of $f( f^{-1}({\bar V'})$ contained  inside $V'''.$ Hence 
       for $h\in {\widetilde W}(f )$ with  $\|h-f\|<\epsilon'$ with $ \epsilon' $ small enough, the image $h(f^{-1}({\bar V'}))\subset V'''$.   We  will  assume that any element $h$ in  ${\widetilde W}(f )$ already has this property.

       % Thus we may  select an open set $D''$ with $f^{-1}({\bar V'})\subset D''$ and $f(D'')\subset V'''$, hence $h(D'')\subset V'''$ by above argument.  
        Then the non-empty compact set $f(S^2 -D''')$ has no intersection with  $ {\bar V'}$  since otherwise, there exists $x\in S^2 -D'''$ such that $f (x)\in    {\bar V'}$ so that $x\in f^{-1} ({\bar V'}) \subset D'''$ contradiction with $x\in S^2-D'''.$ 
     %Now   choose   $V'''$ above such that $V'''$ has no intersection with a   small open neighborhood of $f(S^2 -D'')$. 
    By the argument before we may assume that   for all $h\in {\widetilde W}(f )$, we still have ${\bar V'}\cap h(S^2 -D''')=\varphi.$ In other words, for all $h\in {\widetilde W}(f )$, the inverse images
     $h^{-1}({\bar V'})\subset D'''.$

       Then  by replace $D$ and $V$ by $D''$ and $V''$, the proposition  implies that $T_{\nu}$ can be considered as a smooth section ${\widetilde { W}}(f; D'')\rightarrow {\widetilde  {\cal T}}_{k-1, p}(f; D'')$, hence a smooth map $T_{\nu}:{\widetilde { W}}(f; D'')\rightarrow L_{k-1}^p(D'', f_{D''}^*TV'')=({\widetilde  {\cal T}}_{k-1, p}(f; D''))_f.$ Let $P_{\beta}: {\widetilde  W}(f )\rightarrow L_{k-1}^p(S^2, {\bf R})$ defined by the pull backs of $\beta$, $P_{\beta}(h)=\beta\circ h.$  We will show in next lemma that
$P_{\beta}$ is of class $C^{\infty}.$ Denote the multiplication  map by  $m:L_{k-1}^p(S^2, {\bf R})\times L_{k-1}^p(D'', f_{D''}^*TV'')\rightarrow L_{k-1}^p(D'', f_{D''}^*TV'')$ given by $m(a, \xi)=R(a)\cdot \xi$. Here $R(a)$ is the restriction of $a$
   to $D'$. Then $m$ is smooth. Now consider   $$F=m\circ (P_{\beta}, T_{\nu})\circ (Id, Res): {\widetilde  W}(f )\rightarrow {\widetilde  W}(f )\times {\widetilde { W}}(f; D'')$$ $$  \rightarrow L_{k-1}^p(S^2, {\bf R}) \times L_{k-1}^p(D'', f_{D''}^*TV'')\rightarrow   L_{k-1}^p(D'', f_{D''}^*TV'').$$  Here $Res:{\widetilde  W}(f )\rightarrow {\widetilde { W}}(f; D'') $ is the restriction map. Then  it is smooth. It is easy to check  that $F(h)= (\beta\circ h)\cdot T_{\nu}(h)$ with obvious interpretations of the corresponding domains and ranges. Hence $F$ is essentially
  equal to $ E_{\beta} (T_{\nu}) : {\widetilde  W}(f )\rightarrow L_{k-1}^p(S^2, f^*TM)$ except the range of the map is $L_{k-1}^p(D'', f_{D"}^*TV'').$ However since the support of $\beta$ is contained in $V'\subset \subset V'''$,  the support of $\beta\circ h$ is contained in the compact set $h^{-1}({\bar V'})\subset  D'''.$    Now let $\gamma$ be a cut-off function defined on $S^2$ with support contained in $D''$ and $\gamma=1$ on $D'''$. Then it induces a map $E_{\gamma}:L_{k-1}^p(D'', f_{D''}^*TM)\rightarrow L_{k-1}^p(S^2, f^*TM)$ by multiplying with $\gamma$. Clearly $E_{\gamma}$ is linear
   and continuous,  and  hence smooth. Now $ E_{\beta} (T_{\nu}) : {\widetilde  W}(f )\rightarrow L_{k-1}^p(S^2, f^*TM)$ is equal to $E_{\gamma}\circ F$, hence   smooth as well.
        
        \QED 
        
 \begin{lemma}
 	Let $\beta:M\rightarrow {\bf R}$ be a smooth function. Then the  map $P=P_{\beta}: {\widetilde W}(f) \rightarrow L_k^p(S^2, {\bf R}) $ defined by $P(h)=\beta \circ h$  is of class $C^{\infty}$.
 \end{lemma}       
 
 \proof       
        
  The proof is a simple application the theorem on smoothness of section functor in [P]. The two relevant  finite dimensional bundles are 
 $p_1:B_1\subset E_1=f^*TM\rightarrow S^2$ and   the trivial bundle  $E_2=S^2\times {\bf R}^1\rightarrow S^2$.
   The bundle map  $p=p_{\beta }:B_1\rightarrow E_2$ is defined by $p(\xi)=(x,  \beta (exp_{f(x)}\xi))$ for $\xi\in B_1$ with $p_1(\xi)=x\in S^2.$  
    
    Then by    the  theorem on section functor, ${\Gamma}_{k, p}(p_{\beta }):{\Gamma}_{k, p} (B_1)= {\widetilde W}(f)\rightarrow {\Gamma}_{k, p} (E_2)\simeq L_{k}^p(S^2, {\bf R})$ is of class $C^{\infty}.$ It is easy to check that $P_{\beta}={\Gamma}_{k, p}(p_{\beta })$, hence of class $C^{\infty}$ as well.
    
     \QED   
         
     \begin{cor}
     	Let $\beta:M\rightarrow {\bf R}^m$ be a smooth function. Then the  map $P=P_{\beta}: {\widetilde W}(f) \rightarrow L_k^p(S^2, {\bf R}^m) $ defined by $P(h)=\beta \circ h$  is of class $C^{\infty}$. 
     \end{cor}

       In particular,
        let  $X $ be  a smooth vector field with support in a local chart $V''$ of $M$.  Assume that $\dim (M) =m$ and ${\bf t}=\{t_1, \cdot, t_m\}$   is a smooth local  frame of $TM$ on $U$.  Then $X=\Sigma_{\nu} X^{\nu}t_{\nu}$ with 
        $X^{\nu}: M\rightarrow {\bf R}$ supported in $V''.$
         Then ${\bf P}=(P_{X^1}, \cdots, P_{X^m}):{\widetilde W}(f)\rightarrow L_k^p(S^2, {\bf R}^m) $  is of class $C^{\infty}.$ Now take
         $X=\beta t_{\nu}$ before. Then the  map ${\bf P}$  here is just the corresponding  map for the section $E_{\beta}(T)$  in the previous
          theorem.  This  proved the section 
        $E_{\beta}(T)$ itself  is smooth except that  we are not using the
         standard local trivialization.

     \medskip
      In the rest of this section, we will use  several basic results whose proofs are much easier for the space $L_k^p(S^2, {\bf R}^m)$ than for the general mapping space ${\cal M}_{k, p}.$ The key step to reduce the proofs of theses results to the case  $L_k^p(S^2, {\bf R}^m)$ is  the following proposition.
      
      \begin{pro}
      Let $\iota: M\rightarrow {\bf R}^m$ be an isometric embedding. Then the induced map $\iota_{*}:	 {\cal M}_{k, p}\rightarrow L_k^p(S^2, {\bf R}^m)$ is a closed (and splitting) $C^{\infty}$ embedding.
     
      \end{pro}
      
      \proof
      
      Unless it is a Hilbert space, for a general $L_k^p$-space, a close subspace may not have a complement. On the other hand, the usual definition of a closed submanifold of a Banach manifold requires the local splitting property so that most of the usual properties of submanifolds in finite dimensional  case can be established accordingly for the infinite dimensional  case.
       For our case here, the following implies the local splitting of the embedding.
      
      For any smooth $f:S^2\rightarrow M,$  $$L_k^p(S^2, {\bf R}^m)\simeq L_k^p(S^2, f^*(T{\bf R}^m))= L_k^p(S^2, f^*(TM))\oplus  L_k^p(S^2, f^*(N_M)).$$ Here $N_M$ is the normal bundle of $M$ in ${\bf R}^m$ defined by 
       $(N_M)_m=\{v\in T_m{\bf R}^m\, |\, <v, T_mM>=0\}$ for any $m\in M$  so that $T_m{\bf R}^m=T_mM\oplus (N_M)_m$.
      
      The rest of the proof is a routine verification. 
      
      \QED

      In  Section 4 we will give a general construction that implies the closeness of the embedding in above  proposition.

      The reparametrization group $G={\bf SL}(2, {\bf C})$ acts continuously on ${\cal  M}ap$. The following 
    are  proved in [L ?] :
    
    (I)  The action of $G$ on  ${\cal  M}ap={\cal  M}ap_{k, p}$
     is 
    $G$-Hausdoroff so that the quotient space ${\cal  M}ap/G$ is  always Hausdoroff. 
    
    \medskip
    (II)  
   A  $L_k^p$-map $f\in {\cal  M}ap$  is said to be  weakly stable with respect to $G$ if its stabilizer $\Gamma_f$ is  compact. Let ${\widetilde{\cal B}}^w$  be the collection of all weakly stable maps in  ${\cal  M}ap$. Then any element other than the constant ones in ${\cal  M}ap$ is always weakly stable. 
   Moreover  the $G$- action on ${\widetilde{\cal B}}^w$  is proper so that for any $f\in {\widetilde{\cal B}}^w$ there is an open neighborhood $U$ of $f$ and a compact subset $K$ of $G$  such that  for any $g\in G\setminus K$ and $h\in U$, $g\cdot h\not \in U.$
    
    \medskip
    
    (III) A weakly stable $L_k^p$-map is said to be stable if its stabilizer is finite. Let ${\widetilde {\cal B}}=:{\widetilde {\cal B}}_{k, p}$ of  the collection of  the stable  $L_k^p$-maps. Then for any $f\in {\widetilde {\cal B}}$, there is a local slice $S_f$  that is transversal to the all   $G$-orbit ${\cal O}_h$ at finitely many (uniformly bounded) points for  
    $h$ sufficiently close to $f$. 
    
    To simplifying our presentation, we  assume that there is  at least one  point at which $f$ is a
    local embedding of class $m_0\geq 1$. This  condition implies that $f$  not only  is stable but also  has a local slice $S_f$ constructed below   using evaluation maps.
    %throughout this paper, we  assume that for any stable map $f$  (1) the stabilizer $\Gamma_f$ is always trivial;  (2)  From  first condition,  the $G$-action on ${\widetilde {\cal B}}$  is free so that the quotient space ${\cal B}={\widetilde {\cal B}}/G$  of the unparametrized stable maps is  a topological Banach manifold with thees local slices   $S_f$  as its  Banach coordinate system. 
     % The second condition implies that $f$  not only  is stable but also  has a local slice $S_f$ constructed  geometrically described below.

       %Note that from above discussion,  we only need to construct a local slice $S_f$ near $f$ for the local   $G$-action on a neighborhood ${\widetilde W}(f)$. 

   %  Its restriction on the local slices can be considered as a topological Banach bundle ${\cal L}={\cal L}_{k-1, p}\rightarrow {\cal B}_{k, p}$ on ${\cal B}.$

    ${\bullet }$ Local slices $S_f$ by evaluation map.
    % Natural local  coordinate charts of ${\widetilde {\cal B}}$  and  

      Recall the definition of the local slice  $S_f$ by evaluation map  for $f\in {\widetilde {\cal B}}_{k, p}$ as  follows.
       % natural coordinate chart ${\widetilde W}(f)$ and

       It is sufficient to assume that $f$ is of class $C^{\infty}.$  Let ${\hat W}= {\hat W}_{\epsilon}$ be the
       $\epsilon$-ball of $L_k^p(S^2, f^*(TM)).$ Then we define ${\widetilde W}_f=:Exp_f {\hat W}$. Here the exponential map is taken with respect to an 
       $f$ dependent metric specified below. Since there is at least one point, and hence any point in a neighborhood of that point,  where
       $f$ is a local embedding, we can fix three of such points as three standard
        marked points ${\bf x}=\{x_1, x_2, x_3\}=\{0, 1, \infty\}$ on $S^2$.
       
       Assume that the metric used to define $Exp_f$ depending on $f$ in the sense that it  is flat  near  $f(x_i)\in M, i=1, 2, 3, $ such  that an Euclidean neighborhood
       $U_{f(x_i)}$ of $f(x_i)$ is identified with a small ball $B_{f(x_i)}\subset  T_{x_i}M$ by $exp_{f(x_i)}:B_{f(x_i)}\rightarrow U_{f(x_i)}$.  By this identification $U_{f(x_i)}$ is   decomposed as  $U_{f(x_i)}=V(f(x_i))\oplus H(f(x_i))$  with two flat summands.  
       Here $V(f(x_i))$ is the local image of $f$ near $x_i$ and $H(f(x_i))$ is the local 
       flat hypersurface of codimension $2$ transversal to $f(x_i)$ at its origin.  To define the local slices, we use the following lemma  in [L]. 
       
       \begin{pro}
       	The  $3$-fold evaluation map at ${\bf x}$,  $ev_{\bf x}:{\widetilde W}_f\rightarrow M^3$ defined by $ev_{\bf x}(h)=(h(x_1), h(x_2), h(x_3))$
       	is a smooth submersion.
       	\end{pro}
       	\proof
       		
       	It is sufficient to look at the case  for the evaluation map $ev_x$ with $x$ being one of the $x_i, i=1, 2, 3.$
       		
       	 We use  the  local charts $Exp_f:{\hat W}\subset L_k^p(S^2, f^*(TM))\rightarrow {\widetilde W}_f$ and $exp_{f(x)}:T_{f(x)}M\rightarrow M$.	 Then under these coordinate charts, $ev_x$ is given by
       	 $${\hat ev}_x:\xi\rightarrow Exp_f\xi\rightarrow  exp_{f(x)}(\xi(x))\rightarrow \xi(x).$$ In other words, 
       	 ${\hat ev}_x$ is just the restriction to ${\hat W}$ of the evaluation map 
       	${\hat ev}_x: L_k^p(S^2, f^*(TM))\rightarrow T_{f(x)}M$	given by $\xi\rightarrow \xi(x)$,   which is linear and continuous,  hence smooth.
       	
       \QED

       	 Denote  the three local flat hypersurfaces $H(f(x_i)), i=1,2,3$ together by ${\bf H}$ and the corresponding flat local images $V(f(x_i)), i=1,2, 3$ by ${\bf V}={\bf V}(f({\bf x}))$ with coordinates ${\bf v}$. Let $\pi_{{\bf V}}=\oplus_{i=1}^3\pi_{V(f(x_i))}$ and 
       	 $\pi_{V_i(f)}:U_{f(x_i)}=V(f(x_i))\oplus H(f(x_i))\rightarrow V(f(x_i))$ be  the projection map. 
       	 
       	 For $\epsilon $ small enough, we define $S_f =(ev_{\bf x})^{-1}({\bf H})=(\pi_{{\bf V}}\circ ev_{\bf x})^{-1}(0). $  Then $S_f$ is a $C^{\infty}$ submanifold of codimension six in 
       	${\widetilde W}(f)$, which is also denoted by 	${\widetilde W}(f;{\bf H}).$
       	
       	Let ${\hat S}_f$ be the "lifting" of $S_f$ in  $ L_k^p(S^2, f^*(TM)$ so that 
       	$S_f=Exp_f ({\hat S}_f). $ Then ${\hat S}_f$  is not an open  ball in  a linear subspace  of $ L_k^p(S^2, f^*(TM)$ in general if
       	the map $Exp_f$ is defined by the fixed $g_J$-metric on $M.$ However,  with respect to the $f$-dependent metric above,   $H (f(x_i))$ is a geodesic submanifold so that for any $\xi \in {\hat W}\subset L_k^p(S^2, f^*(TM))$, $h=Exp_f\xi$ is in $S_f$ if and only if $\xi(x_i)\in T_{f(x_i)}{ H}_i, i=1, 2, 3$. This implies  that ${\hat S}_f$ is the $\epsilon$-ball of the Banach space 
       	$L_k^p(S^2, f^*(TM);{T_{f({\bf x})} {\bf  H}})$. Here $L_k^p(S^2, f^*(TM);{T_{f({\bf x})} {\bf  H}})$ is the linear subspace  consisting of $\xi\in  L_k^p(S^2, f^*(TM))$ such that $\xi(x_i)$   is in the tangent space  $T_{f(x_i)}{ H}_i, i=1, 2, 3$.

       To define the coordinate transformations between local slices, we need the following proposition proved in [L  ] and [L ?].	
       	
       	\begin{pro}
       		The composition of the action map  with the $3$-fold evaluation map, $ev\circ \Psi_{{\widetilde{\cal B}}}:G\times {\widetilde{\cal B}}\rightarrow M^3$,  given by $(g, h)=(h\circ g (x_1),h\circ g (x_2), h\circ g (x_3))$,  is of class $C^{m_0}.$

       	\end{pro}
       	
       Now assume that $f': S^2\rightarrow M$ a stable $L_k^p$-map that is equivalent to $f$ in the sense that (i) there is a boholomorphic map $\phi:S^2\rightarrow S^2$ such that $f'=f\circ \phi;$ (ii) $f'$ is also a local embedding at the standard three marked points  ${\bf x}$. Note that (ii) implies that $f$ is a local embedding at both marked point sets ${\bf x}$ and ${\bf y}=:\phi ({\bf x}).$  Under theses conditions, let $S_{f'}$ be the local slice in 	${\widetilde W}_{f'}$, and $\Psi_S:S_{f}\rightarrow S_{f'}$ be  the coordinate
        transformation. The following is proved in [L].
        
        \begin{pro}
       Assume that $m_0>1.$ Then there is a $C^{m_0}$-smooth function $T:S_f\rightarrow G$ with $T(f)=\phi$ such that for any $h\in S_f$, $\Psi_S (h)=h\circ T(h).$
       
      % Moreover, for $L_1^p$-norm with $p>2$ so that  $m_0=[1-2/p]=0$ above conclusion still hold if the center $f$ is chosen be to of class $C^{\infty}$ (or  $C^l$ with $l>1$ ). In other words in this case $T$ is only continuous.
       
    \end{pro}
    
    \proof
    
    %We only prove the case for $m_0>1$ and leave the case for $L_1^p$-norm to the readers.
    
    Consider the  map  $F=:\pi_{{\bf H}_{f'}^{\perp}}\circ ev_{f'(\bf x)}\circ \Psi_{{\widetilde{\cal B} }}:S_f\times G\rightarrow M^3\rightarrow {\bf H}_{f'}^{\perp},$ where $\Psi_{{\widetilde{\cal B} }}$ is the action map restricted to $S_f\times G $. Here ${\bf H}_{f'}^{\perp}$ is the summand in the local decomposition of a flat neighborhood $U_{f'}$ of $M^3$ near $f'({\bf x})$, $U_{f'}={\bf H}_{f'}\oplus {\bf H}_{f'}^{\perp}$. Note that as before here the local flat  metric used is $f'$-dependent. $\pi_{{\bf H}_{f'}^{\perp}}:U_{f'}\subset M^3\rightarrow  {\bf H}_{f'}^{\perp}$ is the projection with respect to above local decomposition.
     Then  for $(h, g)\in S_f\times G$ ,  $h\circ g$ is  in  $S_{f'}$ if and only if  $F(h, g)=0$. In particular, 
      for  any  $h\in S_f$, there is a $g=T(h)\in G$ such that  coordinate transformation $\Psi_S (h)=h\circ T(h) \in S_{f'}$ so that $(h, T(h)) $ solves the equation $F(h, g)=0.$
       Now under the assumption $m_0>1$, the function $F$ is at least of class $C^{m_0}$
        with $F(f, \phi)=0$. Moreover by the construction the partial derivatives along $G$, $\partial_{g} F|_{(f, \phi)}:T_{\phi }G:\rightarrow (T{\bf H}_{f'(0)}^{\perp})_0\simeq {\bf H}_{f'}^{\perp}$ is surjective. Then by implicity function theorem, there is an unique such $T:S_f\rightarrow G$ of class $C^{m_0}$ with the desired property.

        \QED

        The above argument implies the following corollary.
        
        \begin{cor}
        	Given a local slice $S_f\subset {\widetilde W}(f),$ there is a $C^{m_0}$-smooth function $T:{\widetilde W}(f) \rightarrow G$ such that for any $h\in {\widetilde W}(f)$, $h\circ T(h)\in S_f.$
        \end{cor}
       
       Now let  ${\widetilde W}(f)=\cup_{g\in G_e}g(S_f)$  be the decomposition of  ${\widetilde W}(f)$ by the images of $S_f$ under the local $G$-action. Here
       $G_e$ is the corresponding local group. When $f$ is smooth, we may identify
       the ${\cal O}^{G_e}_f $ with $G_e$ smoothly. Then the above corollary implies
       the following corollary.
       
       \begin{cor}
       There  is a $C^{m_0}$-smooth projection map 
       $ T_f:{\widetilde W}(f)\rightarrow {\cal O}^{G_e}_f\simeq G_e$ such that for any $g\cdot f \in {\cal O}^{G_e}_f $, the fiber $T_f^{-1}(g\cdot f )=T^{-1}(g)$ is $g(S_f)$ for any $g\in G_e$.

       \end{cor}

       ${\bullet }$ Local slices $S_f$ by  using $L^2$-metric on $T_fL_k^p(S^2, {\bf R}^m)$
       
       Still assume that  the center $f$ is  of class $C^{\infty} $ so that    the orbit ${\cal O}_f$ is a smooth submanifold of ${\widetilde{\cal B}}_{k, p}\subset L_k^p(S^2, {\bf R}^m).$ Here ${\widetilde{\cal B}}_{k, p}$ is considered as 
       a closed submanifold of $ L_k^p(S^2, {\bf R}^m)$  induced from an embedding 
       $M\subset   {\bf R}^m$. 
       Then the   tangent space $ T_{f}({\cal O}_f) $ is  a $6$-dimensional  linear subspace  of $$T_f{\widetilde{\cal B}}_{k, p}=L_k^p(S^2, f^*TM)\subset T_{f}L_k^p(S^2, {\bf R}^m)=\{f\}\times L_k^p(S^2, {\bf R}^m).$$
       
       Now for any $\xi=(f, \eta)\in \{f\}\times C^{\infty}(S^2, {\bf R}^m)\subset T_{f}L_k^p(S^2, {\bf R}^m)=\{f\}\times L_k^p(S^2, {\bf R}^m),$
       consider the function $I_{\xi}: L_k^p(S^2, {\bf R}^m)\rightarrow {\bf R}$ defined by 
       $I_{\xi}(h)=<\xi, h-f>_2=:\int_{S^2}<\eta,  h-f> dvol_{S^2}.$ In other words, 
       $I_{\xi}(h)$ is defined by the $L^2$-inner product on $T_fL_k^p(S^2, {\bf R}^m)$ of $\xi$  and  the displacement from $f$ to $h$.
       Let $F_{\xi}=I_{\xi} \circ \Psi:G\times L_k^p(S^2, {\bf R}^m)\rightarrow L_k^p(S^2, {\bf R}^m)\rightarrow {\bf R}$, where $\Psi:G\times L_k^p(S^2, {\bf R}^m)\rightarrow L_k^p(S^2, {\bf R}^m)$ is the action map. Then we have the following proposition.

       \begin{pro}
       The function $F_{\xi}:G\times L_k^p(S^2, {\bf R}^m)\rightarrow {\bf R}$ is of class $C^{\infty}.$	
       \end{pro}

       \proof
       
       For $(\phi, h)\in G\times L_k^p(S^2, {\bf R}^m), $
       $$F_{\xi}(\phi, h)=\int_{S^2}<\xi,  h\circ \phi -f> dvol_{S^2}$$ $$ =\int_{S^2}<\xi,  h\circ \phi> dvol_{S^2} -\int_{S^2}<\xi,f> dvol_{S^2}.$$
       
       Thus upto a constant $$F_{\xi}(\phi, h)=\int_{S^2}<\xi,  h\circ \phi> dvol_{S^2} =\int_{S^2}<\xi\circ \phi^{-1},  h> det^{-1} (\phi) dvol_{S^2}$$
       $$= \int_{S^2}<det^{-1} (\phi)\cdot \xi\circ \phi^{-1},  h>  dvol_{S^2}.$$
       
       Clearly $F_{\xi}$ is a smooth function of class $C^{\infty}$. Indeed, 
       Since $\xi$ is fixed and smooth, the function $A_{\xi}: G\rightarrow C^{\infty }(S^2, {\bf R}^m)\subset 
       L_k^p(S^2, {\bf R}^m)$ given by $A_{\xi}(\phi )=det^{-1} (\phi)\cdot \xi\circ \phi^{-1}$ is smooth of class $C^{\infty}$. This together with the smoothness  of the  $L^2$-paring $<-, ->_2$ on $ L_k^p(S^2, {\bf R}^m)$ implies that 
       $F_{\xi}=<-, ->_2\circ (A_{\xi}, Id)$  is of class $C^{\infty}.$
       \QED

       Let ${\hat {\bf e}}=\{{\hat e}^1, \cdots, {\hat e}^1  \}$ be a basis of the Lie algebra $T_eG$, and $ {\bf e}_f=\{e_f^1, \cdots e_f^6\}$ be  the induced  basis of $T_f({\cal O}_f)$ by the infinitesimal action of $G$ at $e$.
        Consider $I_{{\bf e}_f}=(I_{{ e}^1_f}, \cdots, I_{{ e}^6_f}): L_k^p(S^2, {\bf R}^m)\rightarrow {\bf R}^6$ and 
       $F_{{\bf e}_f}=(F_{{ e}^1_f}, \cdots, F_{{ e}^6_f}):G\times L_k^p(S^2, {\bf R}^m)\rightarrow {\bf R}^6.$ Then both of them are  of class $C^{\infty}$. So are their restrictions  to $L_k^p(S^2, M)$ and $G\times L_k^p(S^2, M)$ respectively.
       
       Now consider $I_{{\bf e}_f}:L_k^p(S^2, M)\rightarrow {\bf R}^6.$  Then (1)$I_{{\bf e}_f}(f)={\bf 0}$; (2) The derivative at $f$, $D_f(I_{{\bf e}_f})$ is surjective since the restriction of $D_f(I_{{\bf e}_f})$ to $T_f({\cal O}_f)$ is. This implies that $S_f^{L_2}=:(I_{{\bf e}_f})^{-1}({\bf 0})$ is a
       $C^{\infty}$ submanifold of ${\widetilde W}(f)$ with codimension six. It follows from the construction that	$S_f^{L_2}$ is  a local slice for the local $G$-action on ${\widetilde W}(f).$
       
       If there is no confusion we will still denote $S_f^{L_2}$ by $S_f$.
       
     Now consider   $F_{{\bf e}_f}:G_e\times {\widetilde W}(f)\rightarrow {\bf R}^6.$
      Then (1)  $F_{{\bf e}_f}(e, f)={\bf 0}$; (2) the partial derivative along $G$-direction at $(f, e)$ is surjective. Here $G_e$ is the local "group" near the identity $e$ (= a small neighborhood of $e$ in $G$). The it follows from the implicit function theorem that there is a  $C^{\infty}$ function  $T: {\widetilde W}(f)\rightarrow G$ such that the submanifold $(F_{{\bf e}_f})^{-1}({\bf 0})$ has the form $\{(T(h), h)|\, h\in {\widetilde W}(f) \}$.
   In other words, for any $h\in {\widetilde W}(f), $ $I_{{\bf e}_f}(h\circ T(h))=F_{{\bf e}_f}(T(h), h)= {\bf 0},$ or equivalently $h\circ T(h)\in S_f$. Moreover, for given $h$ if $h\circ \phi \in S_f$ for some $\phi \in G_e$, then 
  $\phi=T(h).$ This proves the following proposition.

      \begin{pro}
      	Given a local slice $S^{L^2}_f\subset {\widetilde W}(f)$, there is a $C^{\infty }$-smooth function $T:{\widetilde W}(f) \rightarrow G$ such that for any $h\in {\widetilde W}(f)$, $h\circ T(h)\in S_f.$

      \end{pro}
      
     {\bf Note} : There is an obvious more natural way to defined the local slice
      $S^{L^2}_f\subset {\widetilde W}(f)$ by setting $S^{L^2}_f=Exp_f {\hat S}$ where $ {\hat S}\subset L_k^p(S^2, f^*TM) $ is  the $\epsilon$-ball of the $L^2$-orthogonal complement  to  the tangent space $T_f ({\cal O}^{G}_f)\subset L_k^p(S^2, f^*TM)$. Then the above proposition is still true.  However its proof
       is slightly harder than the one above.  This is the reason that we use the above definition here. A proof of the  corresponding proposition as above  for this new slice will be given in the forth coming paper on the $C^1$-smoothness of the equivariant extension in the first construction in [L].

      \begin{cor}
      	There  is a $C^{\infty}$-smooth projection map 
      	$ T_f:{\widetilde W}(f)\rightarrow {\cal O}^{G_e}_f\simeq G_e$ such that for any $g\cdot f \in {\cal O}^{G_e}_f $, the fiber $T_f^{-1}(g\cdot f )=T^{-1}(g)$ is $g(S^{L^2}_f)$ for any $g\in G_e$.

      \end{cor}

      Given $f$ and $f'$ in ${\cal B}$, assume that the intersection  of the  orbits  of the two slices, 
       ${\cal O}_{S_f}\cap {\cal O}_{S_{f'}}\not =\varphi$ so that the transformation between these two  slices,  still denoted 
       by $\Psi_S$ is defined on $S_f\cap \Psi_S^{-1}({S_{f'}})$.

       \begin{cor}
       	Assume that  $k-2/p>0$. Then there is a $C^{\infty}$-smooth function $T:S_f\cap \Psi_S^{-1}({S_{f'}})\rightarrow G$ with $T(f)=\phi$ such that for any $h\in S_f\cap \Psi_S^{-1}({S_{f'}})$, $\Psi_S (h)=h\circ T(h).$

       \end{cor}

       \section{Proof of  the Main Theorem}
       
        In this section  we prove the following two versions of the main theorem.
       
       \medskip
       \noindent
       $\bullet$ First version of the main theorem
       
       Now start with 
       a point-section $$\eta\in C^{\infty}(S^2, f^*TM)\subset ({\widetilde {\cal T}}_{k-1, p})_f=L_{k-1}^p(S^2, f^*TM)$$
       with   $$\eta= \Sigma_{i}\alpha_i\cdot \eta=
       \Sigma_{i, \nu}a_i^{\nu}\cdot t_{i\nu}\circ f  $$ $$ =\Sigma_{i, \nu}a_i^{\nu}\cdot \beta_i\circ f\cdot  t_{i\nu}\circ f=\Sigma_{i, \nu}a_i^{\nu}\cdot E_{\beta_i}(  t_{i\nu})\circ f.   
       $$
       
       Here $a_i^{\nu}\in C^{\infty}(S^2, {\bf R})$ supported in $D_i$ such that $\alpha_i\cdot \eta=\Sigma_{ \nu}a_i^{\nu}\cdot t_{i\nu}\circ f$.

       Recall  the assumption that the local $G$-orbit ${\cal O }_{S_f}$ of the local slice $S_f$  in  ${\widetilde W} (f) $ is ${\widetilde W} (f) $ itself. 
       
       Denote by $\eta_{{\cal O }_{S_f}}$ the local $G_e$-equivariant extension of $\eta$ over  ${\cal O }_{S_f}={\widetilde W} (f).$  We require that  the extension $\eta_{{\cal O }_{S_f}}$ to be defined  satisfies the property that
       $$\eta_{{\cal O }_{S_f}}=\Sigma_{i} (\alpha_i\cdot \eta)_{{\cal O }_{S_f}}=\Sigma_{i, \nu}(a_i^{\nu})_{{\cal O }_{S_f}}\cdot (E_{\beta_i}(  t_{i\nu}))_{{\cal O }_{S_f}}.$$  Thus we only need to defined the equivariant extensions $(a_i^{\nu})_{{\cal O }_{S_f}}$ and 
       $(E_{\beta_i}(  t_{i\nu}))_{{\cal O }_{S_f}}$.

       Since $E_{\beta_i}(  T_{i\nu})$ is already smooth  and $G$-equivariant  defined on ${\cal O }_{S_f}={\widetilde W} (f) $,  
        $(E_{\beta_i}(  t_{i\nu}))_{{\cal O }_{S_f}}$ is simply defined to be $E_{\beta_i}(  T_{i\nu})$.
       
       To define  $(a_i^{\nu})_{{\cal O }_{S_f}}$, let ${\cal R}\rightarrow 
       {\widetilde W} (f)$ be the trivial bundle ${\widetilde W} (f)\times L_{k-1}^p(S^2, {\bf R})\rightarrow {\widetilde W} (f).$ Then 
       $a_i^{\nu}\in C^{\infty }(S^2, {\bf R})\subset L_{k-1}^p(S^2, {\bf R})$
       gives rise a constant section, denoted by $(a_i^{\nu})_{{\widetilde W} (f)}:{\widetilde W} (f)\rightarrow {\cal R}.$ Let $(a_i^{\nu})_{S_f}:S_f\rightarrow {\cal R}|_{S_f}$  be the restriction 
       of $(a_i^{\nu})_{{\widetilde W} (f)}$ and $[(a_i^{\nu})_{S_f}]:S_f\rightarrow  L_{k-1}^p(S^2, {\bf R})$ be the corresponding map under the trivialization of ${\cal R}|_{S_f}$.
       Note that $(a_i^{\nu})_{S_f}$ is still a constant section.
       The $G$-action on ${\cal R}$ is  given by  to be $\phi \cdot (h, a)=(h\circ \phi, a\circ \phi)$ for $\phi\in G_e$ and $(h, a)\in {\cal R}={\widetilde W}(f)\times L_{k-1}^p(S^2, {\bf R}).$
       
       Now assume that $S_f$ is one of the two slices defined in Sec. 2. Then we have a $C^{\infty}$ or $C^{m_0}$-smooth map $T:{\widetilde W} (f)\rightarrow G$ such that for any $h\in {\widetilde W} (f)$, $T(h)\circ h \in S_f$.  The $G$-equivariant  extension  $(a_i^{\nu})_{{\cal O }_{S_f}}$ is defined to be $$(a_i^{\nu})_{{\cal O }_{S_f}}(h)=(T(h))^*((a_i^{\nu})_{S_f}(T(h)\circ h)).$$  
       
       Let $[(a_i^{\nu})_{{\cal O }_{S_f}}]:{\widetilde W} (f)\rightarrow L_{k-1}^p(S^2, {\bf R})$ be the corresponding map under the trivialization of ${\cal R}|_{S_f}$. Then $[(a_i^{\nu})_{{\cal O }_{S_f}}](h) =(T(h))^*(a_i^{\nu})=a_i^{\nu}\circ T(h).$ 
       
       \begin{lemma}
       	Assume that $T:	{\widetilde W} (f)\rightarrow G$  is of class $C^{\infty}$ (or $C^{m_0}$ ). The map $[(a_i^{\nu})_{{\cal O }_{S_f}}]$ is of class $C^{\infty}$ (or $C^{m_0}$ ) so that  the $G$-equivariant extension $(a_i^{\nu})_{{\cal O }_{S_f}}$  is  a   smooth section.
       \end{lemma}
       
       \proof
       
       Note that the action map
       
       $\Psi_{(-l)}:G\times  L_{k+l}^p(S^2, {\bf R})\rightarrow  L_{k}^p(S^2, {\bf R})$ is of   class  $C^l$. 
       This implies that for a fixed  $a\in C^{\infty}(S^2, {\bf R})\subset L_{k+l}^p(S^2, {\bf R})$, the map $G\rightarrow L_{k}^p(S^2, {\bf R})$
       given by $\phi\rightarrow a\circ h$ is of class $C^{l}$ for any $l$, hence of class $C^{\infty}.$ 
       Now for $a\in C^{\infty}(S^2, {\bf R})$,  $[a_{{\cal O }_{S_f}}]$ is  the composition of  two smooth maps, $h\rightarrow T(h)$ from ${\widetilde W} (f)$ to $G$
       and $\phi\rightarrow a\circ \phi$ from $G$ to $L_{k}^p(S^2, {\bf R})$.
       Hence $[a_{{\cal O }_{S_f}}]$ is smooth.
       
       \QED
       
        This proves the first version of the main theorem.
        
        \begin{theorem}
        	Given a smooth stable map $f$ and a smooth section 
        	
        	$\eta\in L_{k-1}^p(S^2, f^*TM)$, its  $G$-equivariant extension $\eta_{{\cal O}_{S_f}}$ above as  a section of 
        	${\widetilde {\cal T}}_{k-1, p}\rightarrow {\widetilde W} (f)$ has  the same degree of smoothness as $T$ has. Hence it is either of class $C^{m_0}$ or of class $C^{\infty}$ accordingly.
        \end{theorem}
        
     \medskip
     \noindent
     $\bullet$ Second  version of the main theorem   
       
       Similarly, for  a given point-section $$\xi\in C^{\infty}(S^2, \wedge^{0, 1}(f^*TM))\subset ({\widetilde {\cal L}}_{k-1, p})_f=L_{k-1}^p(S^2, \wedge^{0, 1}(f^*TM))$$
       with   $$\xi= \Sigma_{i}\alpha_i\cdot \xi=
       \Sigma_{i, \nu}{\gamma}_i^{\nu}\cdot t_{i\nu}\circ f  $$ $$ =\Sigma_{i, \nu}{\gamma}_i^{\nu}\cdot \beta_i\circ f\cdot  t_{i\nu}\circ f=\Sigma_{i, \nu}{\gamma}_i^{\nu}\cdot E_{\beta_i}(  t_{i\nu})\circ f.   
       $$
       
       Here ${\gamma}_i^{\nu}\in C^{\infty}(S^2,  \wedge^{0, 1})$ supported in $D_i$ such that $\alpha_i\cdot \xi=\Sigma {\gamma}_i^{\nu}\cdot t_{i\nu}\circ f$.
       
       As  before, the $G$-equivariant extension $$\xi_{{\cal O }_{S_f}}=\Sigma_{i} (\alpha_i\cdot \xi)_{{\cal O }_{S_f}}=\Sigma_{i, \nu}({\gamma}_i^{\nu})_{{\cal O }_{S_f}}\cdot (E_{\beta_i}(  t_{i\nu}))_{{\cal O }_{S_f}}.$$  We  only need to defined the equivariant extensions $({\gamma}_i^{\nu})_{{\cal O }_{S_f}}.$ 
       
       Since the bundle $  \wedge^{0, 1}=:\wedge^{0, 1}_{S^2}\rightarrow S^2$ is not trivial, we introduce an intermediate step as follows.
       
       Note that the group $G$ can be considered as the parametrized moduli space ${\widetilde {\cal M}}(G)=:G$  consisting of holomorphic maps $g:S^2\rightarrow M=S^2$ of class $A=[id]=1\in H_2(S^2).$ Let ${\bf \omega}^{0, 1}={\widetilde {\cal M}}(G)\times C^{\infty}(S^2,\wedge^{0,1})$ be the 
       trivial bundle over ${\widetilde {\cal M}}(G)$. Here the norm used on the fiber $C^{\infty}(S^2,\wedge^{0,1})$ is the $L_k^p$-norm.  The group acts on $C^{\infty}(S^2,\wedge^{0,1})$ by pull-backs, $(g, \gamma)\rightarrow g^*(\gamma)$. Thus 
       the natural action of  $G$  on ${\widetilde {\cal M}}(G)$  lifts to an  action on the  bundle  ${\bf \omega}^{0, 1}$.  Given an $\gamma \in C^{\infty}(S^2,\wedge^{0,1})=C^{\infty}(S^2,\wedge^{0,1})_e$ as a section over the  base point $e\in G={\widetilde {\cal M}}(G)$, let $\gamma_{G}:{\widetilde {\cal M}}(G)\rightarrow {\bf \omega}^{0, 1}$ be its $G$-equivariant extension and  $[\gamma_{G}]:{\widetilde {\cal M}}(G)\rightarrow C^{\infty}(S^2,\wedge^{0,1})$ be the corresponding map under the trivialization. 
       
     %  Note that the group $G$ can be considered as the parametrized moduli space ${\widetilde {\cal M}}(G)=:G$  consisting of holomorphic maps $g:S^2\rightarrow M=S^2$ of class $A=[id]=1\in H_2(S^2).$ Let ${\bf \omega}^{0, 1}={\widetilde {\cal M}}(G)\times C^{\infty}(S^2,\wedge^{0,1})$ be the 
       %trivial bundle over ${\widetilde {\cal M}}(G)$. Here the norm used on the fiber $C^{\infty}(S^2,\wedge^{0,1})$ is the $L_k^p$-norm or   $C^{m_0}+1$-norm.  The group acts on $C^{\infty}(S^2,\wedge^{0,1})$ by pull-backs, $(g, \gamma)\rightarrow g^*(\gamma)$. Thus 
       %the natural action of  $G$  on ${\widetilde {\cal M}}(G)$  lifts to an  action on the  bundle  ${\bf \omega}^{0, 1}$.  Given an $\gamma \in C^{\infty}(S^2,\wedge^{0,1})=C^{\infty}(S^2,\wedge^{0,1})_e$ as a section over the  base point $e\in G={\widetilde {\cal M}}(G)$, let $\gamma_{G}:{\widetilde {\cal M}}(G)\rightarrow {\bf \omega}^{0, 1}$ be its $G$-equivariant extension and  $[\gamma_{G}]:{\widetilde {\cal M}}(G)\rightarrow C^{\infty}(S^2,\wedge^{0,1})$ be the corresponding map under the trivialization. 

       \begin{lemma}
       	For any $\gamma \in C^{\infty}(S^2,\wedge^{0,1})$, the map $[\gamma_{G}]$ is of class $C^{\infty}.$ Thus the $G$-equivariant extension $\gamma_{G}$ of the point-section $\gamma$ is of class $C^{\infty}.$
       	
       \end{lemma}
       
       \proof

       Let  $\wedge^{0,1}=\wedge^{0,1}_{U_1}\cup\wedge^{0,1}_{U_2}\rightarrow S^2=U_1\cup U_2$ be the usual trivialization of he bundle $\wedge^{0,1}$ with respect to the standard covering of $S^2$.  Using the partition on  unit $\beta_1+\beta_2=1$ with respect to the cover,  $\gamma=\beta_1\gamma+\beta_2\gamma.$ Then $\gamma=\beta_1\gamma+\beta_2\gamma=a_1t_1+a_2t_2.$ Here $t_i, i=1, 2$
       is a  smooth complex frame of $\wedge^{0,1}_{U_i}$ and $a_i$ is a complex valued smooth function supported in $U_i$.  Now apply  the   argument before 
        for the smoothness of $E_{\beta_i}(T_{i\nu})$ to this simpler case by replacing $TM$ and $t_{i\nu}$ by $\wedge^{0,1}$ and $t_i, i=1, 2$. Then  we get the $G$-equivariant extension $[\gamma_{G}]:G\rightarrow C^{\infty}(S^2,\wedge^{0,1})\subset L_{k-1}^p$ is of class $C^\infty.$

       \QED
       
       Now define the trivial bundle  ${\bf  \Omega}^{0, 1}|_{{\cal O}_{S_f}}\rightarrow {{\cal O}_{S_f}}$
       to be  the pull-back of 
       ${\bf \omega}^{0, 1}={\widetilde {\cal M}}(G)\times C^{\infty}(S^2,\wedge^{0,1})\rightarrow {\widetilde {\cal M}}(G)$ by the smooth or $C^{m_0}$-smooth  map $p_{{\cal O} _{S_f}}:{\cal O} _{S_f}\rightarrow  G={\widetilde {\cal M}}(G)$. Here  $p_{{\cal O} _{S_f}}$ is defined by $p_{{\cal O} _{S_f}}(h)
       =T(h)$ for $h\in {\cal O} _{S_f}.$  Thus  for $\gamma_i^{\nu}\in  C^{\infty}(S^2,\wedge^{0,1})$, the $G$-equivariant extension of $(\gamma_i^{\nu})_{S_f}$,  $(\gamma_i^{\nu})_{{\cal O} _{S_f}}$, is the same 
       as the pull-back $p_{{\cal O} _{S_f}}^*((\gamma_i^{\nu})_G)$ of the smooth $G$-equivariant section $(\gamma_i^{\nu})_G$. Here $(\gamma_i^{\nu})_{S_f}$ is defined to be the pull-back  $p_{{S_f}}^*(\gamma_i^{\nu})$. Therefore  $(\gamma_i^{\nu})_{{\cal O} _{S_f}}$ is a smooth or $C^{m_0}$-smooth section.

       This proves  the second version the  main theorem.
       
       \begin{theorem}
       	Given a smooth stable map $f$ and a smooth section 
       	
       	$\xi\in L_{k-1}^p(S^2, \wedge^{0,1}(f^*TM))$, its "geometric" $G$-equivariant extension $\xi_{{\cal O}_{S_f}}$ above has  the same degree of smoothness as $T$ has. Hence it is either of class $C^m_0$ or of class $C^{\infty}$ accordingly.
       \end{theorem}

        \section{ The mapping space  ${\cal M}ap_{k, p}(S^2, M)$ as Euler class  of the section $\Delta$}	
        There are several basic results used in this paper whose proof are much easier
        for the Banach space $L_k^p(S^2, {\bf R}^m)$ than for the space  ${\cal M}ap_{k, p}(S^2, M) =L_k^p(S^2, M)$ of  $L_k^p$-maps from $S^2$ to a  compact symplectic   manifold $M$.  We will introduce a general  construction that realizes ${\cal M}ap_{k, p}(S^2, M) $ 
        %as a closed submanifold of the Banach space $L_k^p(S^2, {\bf R}^m)$ for an isometric embedding of $M$ in ${\bf R}^m$. If fact  ${\cal M}ap_{k, p}(S^2, M) $ will be realized 
        as the zero locus of a smooth section $\Delta$ of the bundle ${\cal N}\rightarrow {\cal T}$ where ${\cal T}$ is an tubular neighborhood of ${\cal M}ap_{k, p}(S^2, M) $ in $L_k^p(S^2, {\bf R}^m)$.

         Let $\phi:M\rightarrow {\bf R}^m$ be an isometric embedding  and $p_M:N_M\simeq (T{\bf R}^m)|_M/TM\rightarrow M$  be the normal bundle of $M$ in ${\bf R}^m$.  Denote the corresponding bundle of $\epsilon$-balls in $N_M$ by $p_{B_M}:B_M\rightarrow M $. When $\epsilon$ small enough, the "exponetial' map $exp_{B_M}: B_M\rightarrow {\bf R}^m$  maps the bundle  $B_M$ to the tubular neighborhood $T_M$  of $M$.
       %Clearly the space $L_k^p(S^2, T_M)$ is an open set in  $L_k^p(S^2, {\bf R}^m)$ so that the results can be proved  for $L_k^p(S^2, T_M)$ in essentially the same way as for the Banach space $L_k^p(S^2, {\bf R}^m)$.  The key step to reduce the proofs for general mapping space ${\cal M}ap_{k, p}(S^2, M) $ to 
        %the space $L_k^p(S^2, T_M)$ is the following construction.
        Let $p_{N_{B_M}}:N_{B_M}\rightarrow B_M$ be the pull-back of the normal bundle $p_{N_M}:N_M\rightarrow M$ by the map $p_{B_M}: B_M\rightarrow M$.  Using the inverse of the identification map $exp_{B_M}:B_M\rightarrow T_M$ to pull-back the bundle $p_{N_{B_M}}:N_{B_M}\rightarrow B_M$, we get the bundle $p_{N_{T_M}}:N_{T_M}\rightarrow T_M$ over the tubular neighborhood $T_M$. The tautological section $\delta_{B_M} : B_M\rightarrow  N_{B_M}$ given by $\delta_{B_M}(b)=(b, b)$ gives rise the corresponding smooth section $\delta_{T_M}:T_M\rightarrow N_{T_M}$  that is transversal to the zero section such that $\delta^{-1}_{T_M}(0)=M$. 
        
        Now consider the open set $L_k^p (S^2, T_M)$ of the Banach space $L_k^p (S^2, {\bf R}^m)$. Define the bundle ${\cal N}_{k, p}\rightarrow L_k^p (S^2, T_M)$ by requiring that $({\cal N}_{k, p})_f=L_k^p(S^2, f^*(N_{T_M}))$ for any $f\in L_k^p (S^2, T_M)$. The 
        usual  process in the standard GW theory implies that ${\cal N}_{k, p}$ is indeed a $C^{\infty}$ smooth bundle. The section $ \delta_{T_M}:T_M\rightarrow N_{T_M}$ induces the corresponding smooth
         section $ \Delta:L_k^p (S^2, T_M)\rightarrow {\cal N}_{k, p}.$
        Then by its definition $f\in L_k^p (S^2, T_M)$ lies inside $L_k^p (S^2, M)$ 
        if and only if $\Delta(f)=0$. 
        In other words, $L_k^p (S^2, M)$ is just the zero locus $\Delta^{-1}(0).$ One can show that the section $\Delta$ is transversal to the zero section so that $L_k^p (S^2, M)$ obtained this way is an closed and splitting submanifold of $L_k^p (S^2, {\bf R}^m)$ (see [La] for the definition of splitting submanifold). Indeed, the exponential map $exp_{B_M}:B_M\rightarrow T_M$ induces a corresponding exponential map $Exp_{B_M}:B_{L_k^p (S^2, M)}\rightarrow {\cal T}\subset L_k^2(S^2, T_M)$ where  $B_{L_k^p (S^2, M)}\rightarrow L_k^p (S^2, M)$ is the $\epsilon_1$ ball bundle of ${\cal N}_{k, p}|_{L_k^p (S^2, M)}\rightarrow L_k^p (S^2, M)$ and ${\cal T}$ is  the corresponding open tubular neighborhood of $ L_k^p (S^2, M)$ in $  L_k^2(S^2, T_M)$. In other words,  the restriction  ${\cal N}_{k, p}|_{L_k^p (S^2, M)}\rightarrow L_k^p (S^2, M)$  is just  the normal bundle of 
        $L_k^p (S^2, M)$ in $L_k^p (S^2, {\bf R}^m)$ so that $\Delta$ can be
        interpreted as the corresponding tautological section defined on the  tubular neighborhood ${\cal T}$ as the finite dimensional case before. 
        In particular $\Delta$ is transversal to the zero section. 
        Thus we have realized the infinite dimensional manifold $L_k^p (S^2, M)$ as the "Euler class' of the section $\Delta.$

         Using  above discussion, we now  give a different analytic set-up for GW-theory as follows.
         
         Note that the metric on ${\bf R}^m$ induces a smooth metric on the tangent bundle $T(T_M)\simeq T_M\times {\bf R}^m\rightarrow T_M$.
         Since $exp_{B_M}:B_M\simeq T_M$, $T_M$ itself can be considered as  an open set   of a vector bundle so that $T(T_M)\simeq T(B_M)= V_{B_M}\oplus H_{B_M}$. Here the vertical bundle is the same as $p_{N_{B_M}}:N_{B_M}\rightarrow B_M$ while the horizontal  bundle $ H_{B_M}$ is the bundle  orthogonal to $ V_{B_M}$. Clearly the derivative of projection map $p_{B_M}: B_M\rightarrow M$  identifies each fiber of the  horizontal  bundle $ H_{B_M}$ with the corresponding  fiber of $TM$ so that $ H_{B_M}$ is  identified with the pull-back $p_{B_M}^*(TM)\simeq  H_{B_M}. $ 
         Now using the  identification $T_M\simeq B_M$, we  rewrite 		
         the above decomposition as $T(T_M)= V_{T_M}\oplus H_{T_M}$ with  the projection map  $\Pi_{H}:T(T_M)\rightarrow  H_{T_M}$.
         Then the almost complex structure $J$ and symplectic form $\omega$ becomes the corresponding fiberwise  complex structure $J_{H}$ and symplectic form $\omega_H$ on the bundle  $ H_{T_M}$. 
         		
         Now the linear operator $d:L_k^p(S^2, {\bf R}^m)\rightarrow L_{k-1}^p(S^2, \wedge^1 ({\bf R}^m))$ is obvious a smooth map. It can be interpreted as a smooth  section of the trivial bundle ${\widetilde { \Omega}}_{k-1, p}=L_k^p(S^2, {\bf R}^m) \times L_{k-1}^p(S^2, \wedge^1 (f^*(T{\bf R}^m))\rightarrow L_k^p(S^2, {\bf R}^m)$ for any given
         $C^{\infty}$ "reference" map $f:S^2\rightarrow {\bf R}^m.$   In particular, we may  take the center $f$ to be a smooth map $f:S^2\rightarrow M$. Then  by taking the restriction to the open set $L_k^p (S^2, T_M)$, we get the corresponding smooth section section $d: L_k^p (S^2, T_M)\rightarrow {\widetilde {\Omega}}_{k-1, p}$. Let $d_{H}=(\Pi_{H})_{*}\circ d:L_k^p (S^2, T_M)\rightarrow {\widetilde {\Omega}}_{k-1, p}$ be the composed smooth section. Here $(\Pi_{H})_{*}:
         {\widetilde {\cal T }}_{k-1, p}\rightarrow {\widetilde {\cal H }}_{k-1, p}\subset {\widetilde {\cal T }}_{k-1, p}$ is the $C^{\infty}$ bundle  endomorphism
         of ${\widetilde {\cal T }}_{k-1, p}$ induced  by  $\Pi_{H}; $  and  
        ${\widetilde {\cal T }}_{k-1, p}=L_k^p(S^2, T_M)\times  L_{k-1}^P(S^2, f^*(T{\bf R}^m))\rightarrow L_k^P(S^2, T_M)$ is the trivial bundle  with the "usual" trivialization and ${\widetilde {\cal H }}_{k-1, p}$ is the corresponding subbundle. The  almost complex structure $J_H$ can be interpreted as an endomorphism of $T(T_M)$ by requiring that $J_H=0$ on the vertical bundle. Then $J_H$ induces a $C^{\infty}$ bundle endomorphism $(J_H)_{*}$ of ${\widetilde {\cal T }}_{k-1, p}$ and hence an 
         endomorphism, still denoted by  $(J_H)_{*}$ of the bundle ${\widetilde {\Omega}}_{k-1, p}\rightarrow  L_k^p (S^2, T_M).$ Finally, the complex structure $i_0$ on $S^2$ can be interpreted as a $C^{\infty}$-endomorphism $i_0^*$ of the bundle ${\widetilde {\Omega }}_{k-1, p}\rightarrow  L_k^p (S^2, T_M).$ 
         
         Now  we defined the lifted
         ${\bar {\partial}}_{J_H}$ section of the bundle ${\widetilde {\cal L}}_{k-1, p}\rightarrow  L_k^p (S^2, T_M)$ by ${\bar {\partial}}_{J_H}(h)=
        \Pi_{H} (dh+(J_H)_{*}(h)\cdot i_0^*(dh))$.  Then ${\bar {\partial}}_{J_H}$
         is a smooth section of class $C^{\infty}.$ Note that the section ${\bar {\partial}}_{J_H}$ is in fact  a  section of the smooth susbundle 
        ${\widetilde {\Omega}}^{0, 1}_{H;k-1, p}\rightarrow  L_k^p (S^2, T_M).$
        Thus in therm of the discussion above,  the moduli space of $J$-holomorphic maps,  ${\widetilde {\cal M}}_J=\{f:S^2\rightarrow M|\, \, {\bar{\partial }}_J(f)=df+J(f)\cdot df\circ i_0=0\}$ is just the intersection $\Delta^{-1}(0)\cap {\bar {\partial}}^{-1}_{J_H}(0).$ 
       One can show that the smooth section
        $(\Delta, {\bar {\partial}}_{J_H}): L_k^p (S^2, T_M)\rightarrow
        {\cal N}_{k, p}\times {\widetilde {\Omega}}^{0, 1}_{H;k-1, p}$ is in fact Fredholm with the index of 
        $D_f(\Delta, {\bar {\partial}}^{-1}_{J_H})$ is
        the same as the index of $D_f{\bar {\partial}}_{J}.$
        Thus we have  realized  the moduli space ${\widetilde {\cal M}}_J$ as 
          a Fredholm intersection of two infinite  dimensional cycles in the flat
          ambient space $L_k^p(S^2, T_M).$  
          
          The whole package  for the analytic foundation of GW and Floer type theories can be carried out in this setting and its generalizations.  
          That will lead to simplifications in many cases as here  we are primarily dealing  with the open set ${\cal T}$ in the flat space $L_k^p(S^2, {\bf R}^m).$ 
          
          As an example, the standard basis ${\bf t}=\{t_1, \cdots, t_m\}$ of $ {\bf R}^m$ can be considered as a global frame of the tangent bundle
          $T {\bf R}^m$. Thus it gives rise the smooth global sections $T=\{T_{\nu}, \nu=1, \cdots, m\}$ of the  "tangent" bundle
          $T ({\cal T})_{k-1, p}\rightarrow {\cal T}_{k,p}$.
          This will simplify considerably 
          the discuss of this paper on the $G$-equivariant extension $\xi_{{\cal O}_{S_f}}$ for a point-section $\xi$ if we use this new setting. 
          The details of the discussions in this section  will be treated in a separate
         paper.
        
      \section{Smoothness of $L_k^p$-norm  and related results}  
       In this section we give a simpler proof of the smoothness of the $L_k^p$-norm.
       The  main  results of this section were proved in [C] and reproduced in [L].   Independent proofs with different settings were given in [CLW].
       
       \begin{theorem}
       Assume that $p=2m$ is a positive  even integer. Let $N_k(\xi)= \Sigma_{i=0}^k\int_{\Sigma} |D^i \xi|^pdvol_{\Sigma}.$ Then    $N_k:L_k^p(\Sigma, {\bf R})\rightarrow {\bf R}$ is of class $C^{\infty}	.$ 
       Here $\Sigma$  is an oriented compact Riemannian manifold.
       \end{theorem}
      
      \proof
      
      Note that $N_k=\Sigma_{i=0}^k N^{(i)}$ where $N^{(i)}:L_k^p\rightarrow {\bf R}$ is given by $N^{(i)}(\xi)=\int_{\Sigma} |D^i \xi|^pdvol_{\Sigma}.$
    Then 	$N^{(i)}=N_0\circ D^i:L_k^p\rightarrow L_{k-i}^p\rightarrow R.$
      Here $D^i:L_k^p\rightarrow L_{k-i}^p$ given by $\xi\rightarrow D^i(\xi)$ is 
     linear and  continuous, hence smooth.  Thus $N_k$ is smooth if and only if 
     $N_0$ is.
     
     Now consider the polarization of $N_0.$
    $P:(L_k^p(\Sigma, {\bf R}))_1\times\cdots (L_k^p(\Sigma, {\bf R}))_p\rightarrow {\bf R}$ given by 
     $P(\xi_1, \cdots \xi_{2m})=\int_{\Sigma} < \xi_1, \xi_2>\cdots <\xi_{2m-1 },\xi_{2m}>dvol_{\Sigma}.$ Here each $(L_k^p(\Sigma, {\bf R}))_j$ is a copy of $L_k^p(\Sigma, {\bf R}).$
     
     The key point is that $P$ is well-defined and continuous. Indeed, by Holder inequality
     for $1/p_1+1/p_2+\cdots, +1/p_{2m}=1/r$ with $p_1=p_2=\cdots=p_{2m}=2m=p$ and $ r=1$, $$|P(\xi_1, \cdots \xi_{2m})|=|\int_{\Sigma} < \xi_1, \xi_2>\cdots <\xi_{2m-1 },\xi_{2m}>dvol_{\Sigma}|$$ $$\leq \int_{\Sigma} | \xi_1|\cdot | \xi_2|\cdots |\xi_{2m-1 }|\cdot |\xi_{2m}|dvol_{\Sigma} \leq
      \|\xi_1||_p \cdot \|\xi_2||_p\cdots \|\xi_{2m}||_p.$$ 
     
       It is well-known in the usual Banach calculus that any continuous multi-linear function like $P$ above is of class $C^{\infty}$ (see Lang's book [La]).
       
       Now $N_0=P\circ \Delta_{2m}:L_k^p\rightarrow (L_k^p(\Sigma, {\bf R}))_1\times\cdots (L_k^p(\Sigma, {\bf R}))_p\rightarrow {\bf R}$, where 
       $\Delta_{2m}:L_k^p\rightarrow (L_k^p(\Sigma, {\bf R}))_1\times\cdots (L_k^p(\Sigma, {\bf R}))_p$ is the diagonal map that is smooth.
       Hence $N_0$ is smooth.
       
       \QED

       \begin{cor}
       	The $l$-th derivative of  $N_0$ is equal to zero for $l>p+1=2m+1.$

       \end{cor}
      
      \proof
       
       The $(p+1)$-th derivative  of the  multi-linear map  $P$ above is equal to zero (see [La]). The diagonal  map $\Delta_{2m}$ is linear so that its first and second derivatives are  a constant map and zero respectively. 
      Since $ N_0=P\circ \Delta_{2m} ,$ the conclusion follows from successively applying the chain rule and product rule for differentiations.
       
       \QED

        Now assume that $\Sigma=S^2.$
       \begin{cor}
       	Let  $\Psi: G\times L_k^p(\Sigma, {\bf R})\rightarrow L_k^p(\Sigma, {\bf R})$ be the action map. Here $G={\bf PSL}(2, {\bf C})$ acting on $\Sigma=S^2$ as the group of reparametrizations. Then  $F_k= N_k\circ \Psi: G\times L_k^p(\Sigma, {\bf R})\rightarrow {\bf R}$ is of class $C^{\infty}$.
       \end{cor}
        
        \proof
       
        We only give the proof for $k=0, 1. $  The  general  case can be proved similarly with more complicated notations.
        
        ${\bullet}$ The case of $k=0.$
        
        $F_0(\phi, \xi)=\int_{\Sigma}  |\xi\circ \phi(x)| ^pdvol_{\Sigma}(x)
        =\int_{\Sigma}  |\xi| ^p det^{-1} (\phi)  dvol_{\Sigma}.$
        
        Consider the "polarization" of $F_0,$
        $P_0:G\times (L^p(\Sigma, {\bf R}))_1\times\cdots (L^p(\Sigma, {\bf R}))_p\rightarrow {\bf R}$ given by 
        $P_0(\phi, \xi_1, \cdots \xi_{2m})=\int_{\Sigma} < \xi_1, \xi_2>\cdots <\xi_{2m-1 },\xi_{2m}>det^{-1} (Jac_{\phi})dvol_{\Sigma}.$ As before, $P_0$  is well-defined. Indeed 
        $$|P_0(\phi, \xi_1, \cdots \xi_{2m})|\leq \|det^{-1} (Jac_{\phi})\|_{C^0}\int_{\Sigma} | \xi_1|\cdot | \xi_2|\cdots |\xi_{2m-1 }|\cdot |\xi_{2m}|dvol_{\Sigma}$$ $$  \leq \|det^{-1} (Jac_{\phi})\|_{C^0}
        \|\xi_1||_p \cdot \|\xi_2||_p\cdots \|\xi_{2m}||_p.$$ 
      
       Note that $P_0$ is a composition of two smooth maps $ G\times (L^p(\Sigma, {\bf R}))_1\times\cdots (L^p(\Sigma, {\bf R}))_p\rightarrow (L^p(\Sigma, {\bf R}))_0\times (L^p(\Sigma, {\bf R}))_1\times\cdots (L^p(\Sigma, {\bf R}))_p\rightarrow {\bf R}$ given by 
       $$(\phi, \xi_1, \cdots \xi_{2m})\rightarrow (det^{-1} (Jac_{\phi}), \xi_1, \cdots \xi_{2m})$$ $$ \rightarrow \int_{\Sigma} < \xi_1, \xi_2>\cdots <\xi_{2m-1 },\xi_{2m}>det^{-1} (Jac_{\phi})dvol_{\Sigma}.$$

        Hence $P_0$ and $F_0$ is of class $C^{\infty}.$
      
       ${\bullet}$ The case of $k=1.$
       
       Note that $F_1=F_0+F^{(1)}$. Here $F_1^{(1)}$ is defined by 
       $$F_1^{(1)}(\phi, \xi)=\int_{\Sigma}  |\nabla (\xi\circ \phi)| ^pdvol_{\Sigma} =\int_{\Sigma}  |(\nabla \xi)\circ \phi(x)\cdot Jac_{\phi} (x) | ^pdvol_{\Sigma}(x)
       $$ $$ = \int_{\Sigma}  |\nabla \xi| ^p     \cdot |Jac_{\phi}\circ  \phi^{-1}|^pdet^{-1} (Jac_{\phi})  dvol_{\Sigma}.$$
       
      Consider the "polarization" of $F^{(1)},$
       $P^{(1)}:G\times (L_1^p(\Sigma, {\bf R}))_1\times\cdots (L_1^p(\Sigma, {\bf R}))_p\rightarrow {\bf R}$ given by 
       $P(\phi, \xi_1, \cdots \xi_{2m})=\int_{\Sigma} < \nabla\xi_1, \nabla\xi_2>\cdots <\nabla\xi_{2m-1 },\nabla\xi_{2m}> \cdot <(Jac_{\phi}\circ  \phi^{-1}), (Jac_{\phi}\circ  \phi^{-1})> ^m   \cdot   det^{-1} (Jac_{\phi})dvol_{\Sigma}.$  Then  $P^{(1)}$ is well-defined.  As  before, it is a composition of two smooth  maps and hence smooth.   This  implies that $F_1$ is smooth.

      \QED

\end{document}